\documentclass[onefignum,onetabnum]{siamart171218}
\usepackage[utf8]{inputenc}
\usepackage{algorithmic}
\usepackage{amsmath,amssymb}
\usepackage{mathscinet}

\newcommand{\mathC}{\mathbb{C}}
\newcommand{\mathR}{\mathbb{R}}
\newcommand{\spec}{\mathrm{spec}}

\newcommand{\cg}{\mathrm{cg}}
\newcommand{\cut}[1]{}

\title{A flexible short recurrence Krylov subspace method for matrices arising in the time integration of port Hamiltonian systems and ODEs/DAEs with a dissipative Hamiltonian} %dominant positive definite Hermitian part}
\author{Malak Diab, Andreas Frommer and Karsten Kahl\thanks{University of Wuppertal, Department of Mathematics, Gauss-Strasse 20, 42097 Wuppertal, Germany (\email{mdiab@uni-wuppertal.de}, \email{frommer@uni-wuppertal.de}, \email{kkahl@uni-wuppertal.de} ).} }
\date{\today}

\begin{document}
\maketitle\thispagestyle{plain}
\markboth{A. FROMMER AND K. KAHL}{A FLEXIBLE SHORT RECURRENCE METHOD FOR MATRICES ARISING IN PORT HAMILTONIAN SYSTEMS}

\begin{abstract} For several classes of mathematical models that yield linear 
systems, the splitting of the matrix into its Hermitian and skew Hermitian 
parts is naturally related to properties of the underlying model. This is particularly so for discretizations of dissipative Hamiltonian ODEs, DAEs and port Hamiltonian systems where, in addition, the Hermitian part is positive definite or semi-definite. It is 
then possible to develop short recurrence optimal Krylov subspace methods in which 
the Hermitian part is used as a preconditioner. In this paper we develop 
new, right preconditioned variants of this approach which as their crucial new feature
allow the systems with the Hermitian part to be solved only approximately 
in each iteration while keeping the short recurrences. This new class of methods is particularly efficient as it allows, for example, to use few steps of a multigrid solver or a (preconditioned) CG method for the Hermitian part in each iteration. We illustrate this with several numerical experiments for large scale systems.
\end{abstract}


\section{Introduction}
Consider the linear system
\[
A x = b, \enspace A\in \mathC^{n \times n},
\]
and the splitting
\[
A = H+ S
\]with $H = \tfrac{1}{2}(A+A^*)$ the Hermitian and $S= \tfrac{1}{2}(A-A^*)$ the skew Hermitian part of $A$, $H^* = H, S^* = -S$. Here, ${}^*$ denotes the adjoint with respect to the standard inner product, i.e., for any matrix $B$ the matrix $B^*$ is obtained by transposing $B$ and taking complex conjugates for each entry.

In several important applications, and in particular in linear systems arising from discretizations of dissipative Hamiltonian ODEs and certain port Hamiltonian systems, the hermitian part $H$ is positive definite and ``dominates'' $S$ in the sense that $H^{-1}S$ or $S$ is small. Then, $H$ is an efficient preconditioner for $A$, and as we will discuss in Section~\ref{basics:sec}, there exist optimal short recurrence Krylov subspace methods using this preconditioner.  

Each iterative step of these methods requires one exact solve with the matrix $H$. The purpose of this paper is to go one step further and investigate practically important short recurrence methods where the solves with $H$ are performed only approximately. For example, we might want to use an incomplete Cholesky factorization in cases where an exact factorization is too expensive, or we might want to perform just some steps of (plain or preconditioned) CG or a multigrid solver for $H$. 
We believe that this is an aspect of primordial importance: As long as systems are small enough such that direct factorizations with $A$ are feasible, there is no real need to consider an iterative method which, after all, requires again a direct factorization of a matrix, $H$, of the same size.  
Thus we expect that variants which allow for inexact solves of systems with $H$ to have practical impact.   

The paper is organized as follows: In section~\ref{basics:sec} we review the existing methods with exact solves for $H$ and relate them to a known method for systems where the matrix is Hermitian plus a purely imaginary multiple of the identity. In section~\ref{flexible:sec} we then develop right preconditioned variants of the methods existing in the literature and use those to develop flexible variants. We describe these algorithmically and analyze some of its most relevant properties. In section~\ref{numerics:sec} we  present results of numerical experiments for four examples demonstrating the efficiency gains of the flexible methods over those relying on exact solves. 

\section{Review of methods with exact solves for $H$}
\label{basics:sec}
From now on we assume that the Hermitian part $H$ of $A$ is positive definite. The method proposed by Concus and Golub \cite{ConGol1975} and, independently, by Widlund \cite{Widlund78} as well as the method of Rapoport \cite{Rapoport78} are, mathematically, equivalent to the full orthogonalization method FOM (see \cite{Saad03}, e.g.) and the minimal residual method MINRES \cite{PaiSau75}, respectively, for the preconditioned system
\[
(I+H^{-1}S) x = H^{-1}b
\]
with the standard inner product replaced by the $H$-inner product $\langle x,y \rangle_H := y^*Hx.$ 

In this section, we present the essentials of both these methods. The starting point is to observe is that the matrix $H^{-1}S$ in the preconditioned matrix $H^{-1}A = I+H^{-1}S$ is $H$-anti-selfadjoint according to the following definition and subsequent proposition (whose proof is straightforward).

\begin{definition} A matrix $B\in \mathC^{n \times n}$ is $H$-selfadjoint if $\langle Bx,y \rangle_H = \langle x, By \rangle_H$ for all $x,y \in \mathC^n$. It is $H$-anti-selfadjoint if $\langle Bx,y \rangle_H = - \langle x, By \rangle_H$ for all $x,y \in \mathC^n$.
\end{definition}

\begin{proposition} A matrix $B\in \mathC^{n \times n}$ is $H$-selfadjoint iff
\[
B^* = HBH^{-1}.
\]
It is $H$-anti-selfadjoint if
\[
B^* = -HBH^{-1}.
\]
\end{proposition}

Because $I+H^{-1}S$ is a shifted $H$-anti-selfadjoint matrix, the Arnoldi process for the $H$-inner product produces a three term recurrence. Introducing the notation $\mathcal{K}_m(A,b)$ for the Krylov subspace $\mathcal{K}_m(A,b)= \text{span}(b,Ab,\ldots,A^{m-1}b)$ of order $m$ generated by the matrix $A$ and the vector $b$, this short recurrence property can be formulated as
\[
(I+H^{-1}S)V_m = V_{m+1}T_{m+1,m},
\]
where the columns $v_i$ of $V_m$ form an $H$-orthonormal basis of $\mathcal{K}_m := \mathcal{K}_m(I+H^{-1}S,v_1)  $ and $T_{m+1,m} \in \mathR^{(m+1)\times m}$ is tridiagonal. 
Denoting $T_{m}$ the matrix obtained from $T_{m+1,m}$ by removing the last row, we also have that $T_{m,m} - I$ is skew-symmetric. The iterates of the method of Concus and Golub / Widlund are now variationally characterized by a Galerkin condition which we denote $\ell$GAL to stress the left preconditioning; the iterates of the method of Rapoport are characterized my a minimal residual condition denoted $\ell$MR:
\begin{align}
& \mbox{$\ell$GAL:}  &&r_m := b-(I+H^{-1}S)x_m \perp_H K_m(I+H^{-1}S,r_0),\label{eq:FOM_def}\\
& \mbox{$\ell$MR:}&&\|r_m\|_H = \min_{x \in x_0 +K_m} \|  b-(I+H^{-1}S)x_m \|_H.\label{eq:GMRES_def}
\end{align}
%
%\mbox{ resp.} 
This gives iterates $x_m = x_0 + V_m \zeta_m$ with
\begin{align}
&\mbox{$\ell$GAL:} &&T_m\zeta_m = \|r_0\|_H e_1, \label{fom:eq} \\
&\mbox{$\ell$MR:} &&\zeta_m = \mbox{argmin} \big\{ \big\| \|r_0\|_H e_1 - T_{m+1,m}\zeta_m \big\|_2\big\}. \label{gmres:eq}
\end{align}
Note that in the formula for the $m$-th iterate of the $\ell$MR variant we find the standard $2$-norm in the least squares problem for $\zeta_m$.

%\begin{remark} 
The short recurrences for computing the iterates of $\ell$GAL and $\ell$MR with the $H$-inner product arise using an $LU$ or a $QR$ factorization of $T_m$ or $T_{m+1,m}$ which is updated from one iteration to the next. All eigenvalues of $T_m$ are of the form $1+i\mu$ with $\mu \in \mathR$ so that $T_m$ as well as all its principal minors cannot become singular, which guarantees that a non-pivoted LU-factorization always exists. For the $\ell$MR variant one typically uses a $QR$ factorization of $T_{m+1,m}$ with an implicit representation of $Q$ as a sequence of Givens rotations and $R$ having three non-zero diagonals (the diagonal and two superdiagonals). Overall this process is analogous to the implementations proposed for SYMMLQ and MINRES in~\cite{PaiSau75}, to the implementation of the GAL, MR and ME methods in \cite{Freund90}, to the implementation of $\ell$GAL in \cite{ConGol1975}, \cite{Widlund78}  and of $\ell$MR in \cite{Rapoport78} and to the implementation of QMR as described  in \cite{FreNac91}. We therefore spare ourselves from reproducing technical details, the principles of which can also be found, e.g., in the textbooks \cite{LieStr13}, \cite{Saad03} and in \cite{SimonciniSzyld07}. 
%\end{remark}

\section{Right preconditioned and flexible methods}
\label{flexible:sec}
If we multiply the original equation $Ax = (H+S)x = b$ with the imaginary unit $i$ we obtain
\begin{equation} \label{eq:basic_lin_sys}
(iH+iS)x = ib.
\end{equation}
In here $iS=:B$ is now Hermitian, and preconditioning with $H^{-1}$ yields a shifted $H$-selfadjoint matrix $H^{-1}(iH+B) = iI + H^{-1}B$. 

In his 1990 paper, Freund \cite{Freund90} considers CG type methods for matrices of the form 
\[
i\sigma I + W, \enspace W^{*} = W, \sigma \in \mathR.
\]
%\cite{Freund90} 
There, three short recurrence Krylov subspace methods, all within the framework of the standard $\ell_2$ inner product are developed: GAL, MR and ME. GAL is variationally characterized by a Galerkin condition, MR by a minimal residual condition, and ME (``minimal error'') minimizes the error in 2-norm over the subspace $x_0+(i\sigma I+W)^*K_m(i\sigma I+W,r_0)$. The ME approach  dates back to \cite{Fri62} for the case $\sigma = 0$ and is less commonly used today. In our numerical experiments it behaved quite similarly as the other methods, but with slightly larger errors and residuals than GAL, so we do not consider this variant any further in this paper.

The methods in \cite{Freund90} can be easily adapted to work with the $H$-inner product rather 
than with the standard inner product. The requirement is that $W$ is then an $H$-selfadjoint 
matrix, which is exactly the situation that arises when we precondition \eqref{eq:basic_lin_sys} 
with $H^{-1}$, yielding $W = H^{-1}B$ (and $\sigma = 1$). For any vector $r$, the 
Krylov subspaces $\mathcal{K}(i \sigma I + H^{-1}B, ir)$ and $\mathcal{K}( \sigma I + H^{-1}S, r)$ are identical. This means that for a given initial guess $x_0$, the defining Galerkin conditions produce identical iterates, independently of whether we consider $(\sigma I + H^{-1}S)x = b $ with initial residual $r = b-(\sigma I + H^{-1}S)x_0$ or $(i  \sigma I + B)x = ib $ with initial residual $ ib-(i\sigma I + B)x_0 = ir$. The same holds for the minimal residual based methods. We emphasize this relation in the following proposition for the case $\sigma = 1$.

\begin{proposition} Freund's GAL method with the $H$-inner product for $(i I + H^{-1}B)x = ib$ is identical to the Concus and Golub /  Widlund method for $(I + H^{-1}S)x = b$ in the sense that it produces the same iterates if we have the same initial guess. Freund's MR with the $H$-inner product for $i I + H^{-1}B$ is, in the same sense, identical to Rapoport's method for $I + H^{-1}S$. 
\end{proposition}
There such equivalent to Freund's ME seems to have been published so far.

\subsection{Right preconditioning}
Right preconditioning of $H+S$ with $H^{-1}$ gives rise to the system 
\[
(I + SH^{-1}) \hat x = b, \enspace x = H^{-1}\hat x.
\]
As with left preconditioning, there is a direct and simple correspondence with the ``Hermitian'' formulation on which we will focus in what follows, i.e.\ %and include the parameter $\sigma$
\[
(i\sigma I + BH^{-1}) \hat x = ib, \mbox{ with } B = iS \mbox{ hermitian} \text{\ and\ } \sigma\in \mathbb{R}.
\]
Herein, $W = BH^{-1}$ is $H^{-1}$-selfadjoint, and we again obtain short recurrence methods GAL, MR, ME by using the $H^{-1}$ inner product in the methods and algorithms of \cite{Freund90}.

\begin{proposition}
There exist short recurrence GAL, MR and ME methods for the right preconditioned system $i\sigma I + BH^{-1}$. These give rise to short recurrence methods for $i\sigma H + B$  when using the $H^{-1}$ inner product.
\end{proposition}

Although it might not be obvious at first sight, the implementation of right preconditioning can be done by investing one multiplication with $B$, one with $H$ and just one solve with $H$ per iteration. 
The idea is to carry an additional vector $\tilde v$ which stores $H^{-1}v$ for the current Lanczos vector $v$. 
Details are given in Algorithm~\ref{lanczos:alg} in which, when solving the linear system \eqref{eq:basic_lin_sys}, we take the initial vector $w$ as the initial residual $ib-(i\sigma H +B)x_0$.  We note that an implementation with similar cost is possible for left preconditioning. 
Also note that in the algorithm we have 
\begin{eqnarray*}
 \langle w, \widetilde{v}_{k-1} \rangle & =&  \langle (i\sigma  I + BH^{-1}) v_k, H^{-1}{v}_{k-1} \rangle \\
 &=& \langle H^{-1} v_k, (-i\sigma  I + BH^{-1}) {v}_{k-1} \rangle \\
 &=& \langle H^{-1} v_k, (i\sigma  I + BH^{-1}) {v}_{k-1} \rangle  \quad \mbox{(since $\langle v_k,v_{k-1} \rangle_{H^{-1}}= 0$) }\\
 &=& \langle H^{-1}v_k, \beta_{k-1}v_k - \alpha_{k-1}v_{k-1} - \beta_{k-2}v_{k-2} \rangle \\
 &=& \beta_{k-1},
 \end{eqnarray*}
%that $\beta_k = \langle w, w\rangle_{H^{-1}}^{1/2}  = \langle v_{k+1}, v_k \rangle_{H^{-1}}$, 
which is why we can use $\beta_{k-1}$ in the $H^{-1}$-orthogonalization of $w$ against $v_{k-1}$ in line~\ref{line:orth}.

\begin{algorithm}
\begin{algorithmic}[1]
\STATE{choose initial vector $w$, compute  $ \widetilde{w} = H^{-1}w, \beta_0 = \langle w, \widetilde{w}\rangle$}
\STATE{put $v_1 = w/\beta_0, \widetilde{v}_1 = \widetilde{w}/\beta_0, v_0 = 0$}
\FOR{$k=1,\ldots,m$}
\STATE{$w = i\sigma v_k + B\widetilde{v}_k$ } \COMMENT{we have $\widetilde{v}_k = H^{-1}v_k$ from the previous iteration} \label{line:mvm} 
\STATE{$\alpha_k = \langle w, \widetilde v_k \rangle$}
\STATE{$w = w - \alpha_k v_k - \beta_{k-1} v_{k-1}$} \label{line:orth}
\STATE{$\widetilde{w} = H^{-1}w$}
\STATE{$\beta_k = \langle w, \widetilde{w} \rangle^{1/2}$}
\STATE{$v_{k+1} = w/\beta_k$}
\STATE{$\widetilde{v}_{k+1} = \widetilde{w}/\beta_{k}$} \COMMENT{now $\widetilde{v}_{k+1} = H^{-1}v_{k+1}$}
\ENDFOR
\end{algorithmic}
\caption{$m$ steps of Lanczos for $i\sigma I+BH^{-1}$ with the $H^{-1}$ inner product. $\langle \cdot,\cdot \rangle$ denotes the standard inner product \label{lanczos:alg}}
\end{algorithm}

One can work with $\sigma I + SH^{-1}$ rather than with 
$i\sigma I + BH^{-1}$ in Algorithm~\ref{lanczos:alg}, adapting 
the matrix vector multiplication in line~\ref{line:mvm} and changing 
the last summand $-\beta_{k-1}v_{k-1}$ to $+\beta_{k-1}v_{k-1}$ in line~\ref{line:orth}. 
The reason for the latter is that due to the fact that $SH^{-1}$ is 
$H^{-1}$-anti selfadjoint we now have $\langle w, \widetilde{v}_{k-1} 
\rangle = - \langle H^{-1}v_k, \beta_{k-1}v_k \rangle$. This observation represents the general transition rule from formulations for the systems with $i\sigma I + BH^{-1}, i\sigma I + H^{-1}B$ or $i\sigma H + B$ to systems with $\sigma I + SH^{-1}, \sigma I + H^{-1}S$ or $\sigma H + S$.
Our presentation will continue to be based on the systems using $B$ rather than $S$.

Right preconditioning is the key to develop flexible variants of the GAL, MR and ME methods. The term ``flexible'' includes ``inexact'' methods in which $\widetilde{w} = H^{-1}w$ is replaced by $\widehat{w} = \widehat{H}^{-1}w$ with a fixed matrix $\widehat{H}$. For example, $\widehat{H}$ may result from an incomplete Cholesky factorization. ``Flexible'' also includes the more general possibility to approximate $H^{-1}w$ via a non-stationary iteration such as (possibly preconditioned) CG which, formally, means that $\widehat{H}$ changes from one approximate solve to the next. 

\subsection{Flexible methods}
Using exact multiplications with $H^{-1}$, the Lanczos algorithm given in Algorithm~\ref{lanczos:alg} for the right preconditioned matrix $i\sigma I + BH^{-1}$ produces the relation
\begin{equation} \label{Arnoldi0:eq}
(i\sigma I + BH^{-1})V_m = V_{m+1}T_{m+1,m},
\end{equation}
with
\[
T_{m+1,m} = \begin{bmatrix} \alpha_1 & \beta_1 &   &   &  \\
\beta_1 & \alpha_2 & \beta_2  & & \\
       & \beta_2   & \ddots  & \ddots  & \\
       &           & \ddots  & \ddots  & \beta_{m-1} \\  
       &           &   & \beta_{m-1} & \alpha_m \\
       &           &      &        & \beta_m 
\end{bmatrix},
\]
where now $V_m$ has $H^{-1}$-orthonormal columns. The transformation of an iterate
\[
\widetilde{x}_m = x_0 + V_m \zeta_m
\]
for the right preconditioned system 
back to the iterate $x_m = H^{-1}\widetilde{x}_m = H^{-1}x_0 + H^{-1}V_m\zeta_m$ of the original system is most easily done via a short recurrence update using the preconditioned Lanczos vectors $\widetilde{v}_m = H^{-1}v_m$
which we compute anyways. The relation \eqref{Arnoldi0:eq} can then be stated as
\begin{equation} \label{Arnoldi1:eq}
(i\sigma H + B)\!\!\!\!\!\underbrace{\widetilde{V}_m}_{=\,H^{-1}V_m}\!\!\!\!\!\! = V_{m+1}T_{m+1,m}.
\end{equation}
We can now divise a flexible method in a manner similar to what has been done for FGMRES, the flexible GMRES method \cite{VuikVorst1994} and, with regard to short recurrences, in flexible QMR \cite{SzyVog01}, e.g. Assume that
$z_k$ is only an approximation to $H^{-1}v_k$ and proceed as follows:
\begin{enumerate}
    \item Compute $w= (i\sigma H+B)z_k$
    \item $H^{-1}$-orthogonalize $w$ against $v_k$ and $v_{k-1}$ {\em approximately}. This means that the inner products $\langle w, v_k \rangle_{H^{-1}}$ and $\langle w, v_{k-1} \rangle_{H^{-1}}$ are obtained as standard inner products with $\widehat w$ where $\widehat w$ is an approximation to $H^{-1}w$. 
    \item Approximately $H^{-1}$-normalize the vector resulting from the approximate $H^{-1}$ orthogonalization.  
\end{enumerate}

Algorithm~\ref{flexlanczos:alg} gives the details for one iteration of the resulting right preconditioned flexible Lanczos process. 
\begin{algorithm}
\begin{algorithmic}[1]
\STATE{choose initial vector $w$, compute  $ \widehat{w} \approx H^{-1}w, \beta_0 = \langle w, \widehat{w}\rangle$}
\STATE{put $v_1 = w/\beta_0, z_1 = \widehat{w}/\beta_0, v_0 = 0$}
\FOR{$k=1,\ldots,m$}
\STATE{$w = (i\sigma H + S) z_k$} \COMMENT{we have $z_k \approx H^{-1}v_k$ from the previous iteration}
\STATE{$\alpha_k = \langle w, z_k \rangle$}
\STATE{$\gamma_k = \langle w, z_{k-1} \rangle$}
\STATE{$w = w - \alpha_k v_k - \gamma_{k} v_{k-1}$}
\STATE{$\widehat{w} \approx H^{-1}w$}
\STATE{$\beta_k = \langle w, \widehat{w} \rangle^{1/2}$}
\STATE{$v_{k+1} = w/\beta_k$}
\STATE{$z_{k+1} = \widehat{w}/\beta_{k}$} \COMMENT{now $z_{k+1} \approx H^{-1}v_{k+1}$}
\ENDFOR
\end{algorithmic}
\caption{$m$ steps of right preconditioned flexible Lanczos for $i\sigma H+B$ with approximate solves for $H$. $\langle \cdot,\cdot \rangle$ denotes the standard inner product \label{flexlanczos:alg}}
\end{algorithm}

In the non-flexible Algorithm~\ref{lanczos:alg} we used $ \langle w, \widetilde{v}_{k-1} \rangle
= \langle H^{-1}v_k,v_k \rangle^{1/2} = \beta_{k-1}$, and we avoided to compute $\langle w, \widetilde{v}_{m-1} \rangle$ explicitly.
In the flexible context the above equality does not hold any more, since $\langle w, \widetilde{v}_{m-1}\rangle$ is computed with the approximation $z_{m-1}$ for $H^{-1}v_{m-1}$, whereas $\beta_m$ is computed with the approximation $\widehat{w}$ for $H^{-1}w$.

With the flexible Lanczos process we have the recurrence relation
\begin{equation} \label{eq:Arnoldi2}
(i\sigma H+B)Z_m = V_{m+1}T_{m+1,m},
\end{equation}
were $T_{m+1,m}$ is again tridiagonal, 
\begin{equation} \label{eq:Tflex}
T_{m+1,m} = \begin{bmatrix} \alpha_1 & \gamma_2 &   &   &  \\
\beta_1 & \alpha_2 & \gamma_3  & & \\
       & \beta_2   & \ddots  & \ddots  & \\
       &           & \ddots  & \ddots  & \gamma_{m-1} \\  
       &           &   & \beta_{m-1} & \alpha_m \\
       &           &      &        & \beta_m 
\end{bmatrix}.
\end{equation}
When used to solve the linear system \eqref{eq:basic_lin_sys}, the initial vector $w$ in Algorithm~\ref{flexlanczos:alg} is the residual $ib-(i\sigma H + B)x_0$, i.e., the first Lanczos vector $v_1$ is
\begin{equation} \label{eq:flex_initial_vector}
v_1 = \frac{1}{\beta_0}r_0 \text{ with } r_0 = ib-(iH + B)x_0, \widehat{w} \approx H^{-1}r_0, \beta_0 = \langle r_0,\widehat{w} \rangle^{1/2}.  
\end{equation}

We emphasize that $V_m$ is now only approximately $H^{-1}$-orthonormal and $Z_m$ is only approximately equal to $H^{-1}V_m$.
We define the flexible GAL and MR iterates in analogy to \eqref{fom:eq} and \eqref{gmres:eq} as $x_m = x_0 + 
Z_m\zeta_m$ with 
\begin{align}
&\mbox{FGAL:} &&T_m\zeta_m = \beta_0 e_1,  \label{eq:fgal} \\
&\mbox{FMR:} &&\zeta_m = \mbox{argmin} \{ \| \beta_0 e_1 - T_{m+1,m}\zeta_m\|_2\}, \label{eq:fmr}
\end{align}
where now $V_m$, $T_{m+1,m}$  and $T_m$ obey  
relation $\eqref{eq:Arnoldi2}$. Computationally, we obtain short recurrences for the iterates $x_m$
by updating factorizations of $T_{m,m}$ and $T_{m+1,m}$ just as in the 
non-flexible case; see the discussion at the end of section~\ref{basics:sec}. Commented Matlab implementations of FMR and FGAL are available from the authors on request.
%as supplementary material for this paper\footnote{\textbf{Link goes here}}. 

Note that the FGAL iterates do not fulfill a strict $H^{-1}$-orthogonality relation for the residuals, nor do the FMR iterates minimize the $H^{-1}$-norm of the residual exactly. But, of course, if we solve for $H$ exactly in the flexible Lanczos process, FGAL and FMR reduce to right preconditioned counterparts of $\ell$FOM and $\ell$MR from \eqref{fom:eq} and \eqref{gmres:eq}, respectively. 
%\begin{remark} 

A flexibly preconditioned CG method (FCG) for Hermitian positive definite matrices using a preconditioner which may change from one iteration to the next has been proposed in 
\cite{GolubYe99} and further analyzed in \cite{Notay2000}. The FGAL method can be used in this situation, too, with $\sigma = 0$. FCG and FGAL are similar in spirit but not identical: FCG constructs search directions by applying the variable preconditioner to the current residual and then $A$-orthogonalizing against a fixed number $s$ of the previous search directions. In contrast, FGAL orthogonalizes the next preconditioned Lanczos vector in the $\ell_2$-inner product against the two previous ones. This is different, even when $s=2$. 
%\end{remark}

\section{Properties} \label{properties:sec}

In this section we first state two known convergence bounds for the left preconditioned methods where systems with $H$ are solved exactly. We then show that a result from \cite{Freund90} can be used to improve one of these bounds, a fact that seems to have gone unnoticed so far. We continue the section considering the flexible variants for which we establish two somewhat more basic but algorithmically important properties.

\begin{theorem} \label{thm:exact_bounds_widlund_rapoport} Let $\lambda \geq 0$ be the smallest value such that the interval $[-i\lambda,i\lambda]$ contains all the (purely imaginary) eigenvalues of the $H$-anti-selfadjoint matrix $H^{-1}S$ and let $x_*$ denote the solution of the system $Ax =b$.  Then
\begin{itemize}
    \item[(i)] For the $\ell$GAL iterates $x_m$ from \eqref{eq:FOM_def}
    we have 
    \begin{eqnarray*}
    \|x_{2m}-x_*\|_H &\leq& \|x_0-x_*\|_H \cdot 2 \left(\frac{\sqrt{1+\lambda^2}-1}{\sqrt{1+\lambda^2}+1}\right)^{m}  \end{eqnarray*}
and 
    \begin{eqnarray*}
    \|x_{2m+1}-x_*\|_H &\leq& \|x_1 - x_*\|_H \cdot 2 \left(\frac{\sqrt{1+\lambda^2}-1}{\sqrt{1+\lambda^2}+1}\right)^{m}.
   \end{eqnarray*}
   \item[(ii)] For the $\ell$MR residuals $r_m = b-Ax_m$ from \eqref{eq:GMRES_def} we have
   \[
      \|b-Ax_m\|_{H^{-1}} \leq \|b-Ax_0\|_{H^{-1}} \cdot 2 \left(\frac{\lambda}{\sqrt{1+\lambda^2}+1}\right)^{m}.  
   \]
\end{itemize} 
\end{theorem}

 Part (i) of the theorem goes back to \cite{Eis83,HagLukYou80,SzyWid93}, and part (ii) to \cite{SzyWid93}; see \cite{GueLieMehSzy21} for a detailed discussion. What seems to have gone unnoticed so far is that a better bound for the $\ell$MR iterates 
results from \cite[Theorem~4]{Freund90}, adapted to $H$-selfadjoint matrices. We first formulate the result for an arbitrary value for the shift $\sigma$ before comparing to  Theorem~\ref{thm:exact_bounds_widlund_rapoport}(ii). Note that the original theorem in \cite{Freund90} 
considers residuals for the matrix $i\sigma I +W$ with $W$ Hermitian and the $\ell_2$-norm. For the residuals belonging to matrices $i\sigma I +W$ with $W$ $H$-selfadjoint it holds in exactly the same manner, but now with the $H$-inner product and associated norm.  When translating back to our original matrix $H + S$, we have $\sigma = 1$, and $H$-norms for the residual w.r.t.\ $i I +H^{-1}B$ turn into $H^{-1}$-norms for the residuals w.r.t.\  $H + S$. This gives us the following theorem.

\begin{theorem} \label{thm:exact_bounds_freund} Let $\spec(H^{-1}S) \subseteq i[\alpha,\beta]$ with $\alpha < \beta$ %, put $a = \frac{2i + \beta+\alpha}{\beta-\alpha}$, 
and let $R$ be the unique solution of
\begin{equation} \label{eq:R}
\frac{1}{2}\left(R+\frac{1}{R}\right) = \frac{\sqrt{\beta^2+1}+\sqrt{\alpha^2+1}}{\beta-\alpha}, \enspace R > 1.
\end{equation}
Then
%\begin{itemize}
%    \item[(i)] Provided $T$ is positive definite For the FOM iterates $x^m$ from \eqref{eq:FOM_def}
%    we have 
%    \[
%    \|x_m-x_*\|_T \leq \|x_0-x^*\|_T \cdot \sqrt{ 1+\frac{\sigma^2(\sqrt{\kappa}-1/\sqrt{\kappa})^2}{4\sigma^2+\alpha^2(\sqrt{\kappa}+1/\sqrt{\kappa})^2} } \cdot \frac{2}{R^k+R^{-k}}.
%   \]
%   \item[(ii)] 
the residuals of the $\ell$MR iterates $x_m$ from \eqref{eq:GMRES_def} satisfy
   \[
      \|b-Ax_m\|_{H^{-1}} \leq \|b-Ax_0\|_{H^{-1}} \cdot \frac{2}{R^k+R^{-k}} .
   \]
%\end{itemize} 
\end{theorem}

Compared to Theorem~\ref{thm:exact_bounds_widlund_rapoport} this theorem allows for a more narrow interval to contain the spectrum of $H^{-1}S$ as it does not need to be symmetric with respect to the real axis. Note, however, that for real matrices $\spec(H^{-1}S)$ is always symmetric with respect to the real axis, so that optimal bounds on the spectrum are of the form $i[-\lambda,\lambda]$ as in Theorem~\ref{thm:exact_bounds_widlund_rapoport}. But even with symmetric bounds for the spectrum, the bound of Theorem~\ref{thm:exact_bounds_freund} is always better than that of Theorem~\ref{thm:exact_bounds_widlund_rapoport} as we will explain now. In order to not unnecessarily burden our discussion, we will omit details of somewhat longish but elementary algebraic manipulations in the paragraph to follow.

We first note that the expression $(\sqrt{\beta^2+1}+\sqrt{\alpha^2+1})/{(\beta-\alpha)}$ 
increases monotonically in $\beta$ and decreases monotonically in $\alpha$. Thus, for
$\lambda \geq \max\{\beta,-\alpha\}$, we have
\[
\frac{\sqrt{\beta^2+1}+\sqrt{\alpha^2+1}}{\beta-\alpha} \leq \frac{\sqrt{\lambda^2+1}}{\lambda}.
\]
The solution $R > 1$ of the equation $R+{1}/{R} = c$ with $c \geq  1$ increases monotonically with $c$.  Consequently, when $\lambda \geq \max\{\beta,-\alpha\}$, the solution $R$ of \eqref{eq:R} is larger or equal than the solution $\widehat{R} = \left(\sqrt{1+\lambda^2}+1\right)/\lambda$ 
of $R+ 1/{R} = {\sqrt{\lambda^2+1}}/{\lambda}$. Now, the factor $\left( \lambda/(\sqrt{1+\lambda^2}+1)\right)^m$ in the bound of Theorem~\ref{thm:exact_bounds_widlund_rapoport}(ii) is precisely ${1}/{\widehat{R}^m}$. The
quantity ${1}/({R^m+R^{-m}})$ in the bound of Theorem~\ref{thm:exact_bounds_freund} is 
thus smaller for two reasons: Because of the presence of $R^{-m} > 0$ and because $\widehat{R} > R$
as soon as $\beta \neq -\alpha$, i.e., $\lambda > \beta$ or $\lambda > -\alpha$.

We now turn to the flexible methods. Conceivably, a detailed convergence analysis, depending on the ``degree of inexactness'' of the preconditioner could be based on prior work for flexible CG \cite{Notay2000}, flexible GMRES \cite{Saa93,VuikVorst1994} and, in particular, flexible QMR \cite{SzyVog01}. We anticipate this to be quite involved and also quite technical, so that we do not address this in the context of the present paper. We rather state and prove two basic properties which are essential for the algorithmic aspects. 

The first property is that in typical situations we cannot encounter unlucky breakdowns in the flexible Lanczos process, Algorithm~\ref{flexlanczos:alg}. An unlucky breakdown at step $m$ means that we get $\beta_m = \langle w, \widehat{w} \rangle = 0$ although $w \neq 0$. If there is 
no unlucky breakdown, the FMR iterates all exist up to the step where we obtain the solution, as in that case
the matrix $T_{m+1,m}$ of \eqref{eq:Tflex} has full rank for all relevant $m$. Though, its square part $T_{m,m}$ need not be non-singular, which means that the FGAL iterate does not necessarily exist. However, as is done in the SYMMLQ-algorithm (see \cite{LieStr13,PaiSau75}), the $QR$-factorization of $T_{m,m}$ can still be updated from that of $T_{m-1,m-1}$, which allows to efficiently update the next FGAL iterate from the previous existing one using the $QR$-factorization.  

\begin{proposition} Assume that the inexact solves $\widehat{w} \approx H^{-1}w$ in Algorithm~\ref{flexlanczos:alg} are done in one of the following ways
\begin{itemize}
    \item[(i)] by performing $k$ steps of the CG method with initial guess 0,
     \item[(ii)] by performing $k$ steps of the preconditioned CG method with initial guess 0 and a Hermitian positive definite preconditioner $K$,              
    \item[(iii)] as $\widehat{w} = Gw$ with $G \in \mathC^{n \times n}$ Hermitian and positive definite.
\end{itemize}
Then, if $w \neq 0$, we have $\beta = \langle w, \widehat{w} \rangle > 0$.
\begin{proof}
In case (i), the initial residual for the CG method is $w$, and with the columns of $W \in \mathC^{n \times k}$ being the Lanczos vectors which form an orthonormal basis of the Krylov subspace $\mathcal{K}_k(H,w)$, we know that $\widehat{w} = V_k(V^{*}_{k}HV_k)^{-1}V_k^*w$; see e.g.\ \cite{LieStr13} or \cite{Saad03}. Thus, $\beta_m = \langle w, \widehat{w} \rangle = \langle V_k^*w, (V_k^*HV_k)^{-1}V_k^*w\rangle$, and this quantity is positive as  $V_k^*HV_k$ is Hermitian positive definite and $V_k^*w = \|w\|_2 e_1$, $e_1$ the first canonical unit vector in $\mathC^k$, is non-zero.

Case (ii) follows in a similar manner as (i), replacing $V_k$ by the matrix of $K$-orthonormal Lanczos vectors for $\mathcal{K}_k(K^{-1}H,K^{-1}w)$. 

In case (iii), as $G$ is Hermitian positive definite and $w \neq 0$, we immediately have $\beta = \langle w, \widehat{w} \rangle = \langle w, Gw \rangle > 0$. 
\end{proof}
\end{proposition}

An important example for case (iii) above is when $G=L^{-1}L^{-*}$ with $L$ a sparse lower triangular matrix arising from an incomplete Cholesky factorization of $H$, $H = L^*L-R$. Such $G$ is also a typical preconditioner $K^{-1}$ for case (ii). Another important example for (iii) are Galerkin based multigrid methods for Hermitian positive definite systems with symmetric pre- and post-smoothing and exact coarsest level solves.

The second property we want to stress deals with $\|b-Ax_m\|_{H^{-1}}$ as a measure of the error of iterate $x_m$.  %The $H^{-1}$-norm of the residual should be considered a more insightful quantity than the $\ell_2$-norm of the residual. 
Using $b-Ax_m = A(x_*-x_m)$, we have
\begin{eqnarray*}
\| b-Ax_m \|_{H^{-1}}^2 &=& \langle A(x_*-x_m), A(x_*-x_m) \rangle_{H^{-1}} \\
&=& 
\langle x_*-x_m, x_*-x_m \rangle_{A^*H^{-1}A} \, = \, \| x_*-x_m \|_{A^*H^{-1}A}^2
\end{eqnarray*}
with 
\[
A^*H^{-1}A = (-i \sigma H + B)H^{-1}(i\sigma H + B) = \sigma^2 H +  BH^{-1}B.
\]
Thus, if $H$ dominates $S$, i.e., if $H^{-1}B$ is small, the matrix ${AH^{-1}A}$ is close to $\sigma^2 H$ (and we have equality if $B = 0$). Consequently, when $H$ dominates $B$, up to the trivial factor $\sigma^2$, the norm $\|b-Ax_m \|_{H^{-1}}$ is a good approximation to the $H$-norm of the error. There are multiple reasons why in the case $B=0$ the $H$-norm of the error---and not its $\ell_2$-norm nor the $\ell_2$-norm of the residual, e.g.---should be considered the most adequate measure for the error, in particular when the linear systems arise from a discretization of an underlying continuous, infinite-dimensional equation; see \cite{Arietal13,MalStr15}, e.g. Even if $B \neq 0$, the $H$-norm can still be regarded as the canonical measure for the error.  By the following proposition, an approximate bound for this norm is available at almost no cost for the FMR iterates.

\begin{proposition} \label{prop:bound} Let $v_k, k=1,\ldots,m,$ denote the vectors produced by the flexible Lanczos algorithm, Algorithm~\ref{flexlanczos:alg}, and let $\varrho_m$
be the value of the minimum in the defining equation \eqref{eq:fmr} of the FMR iterate $x_m$, i.e.\
\[
\varrho_m :=   \min_{\zeta \in \mathC^m} \left\{ \| \beta_0 e_1 - T_{m+1,m}\zeta \|_2\right\}.
\]
Then
\begin{equation} \label{eq:FMR_bound}
\| b-Ax_m \|_{H^{-1}} \leq \|V_{m+1}\|_{H^{-1}} \cdot \varrho_m,
\end{equation}
and, if the flexible Lanczos vectors $z_m$ approximate the ``exactly inverted'' vector $H^{-1}v_m$ to a relative accuracy of $\varepsilon$ in the $H$-norm, i.e.\
\begin{equation} \label{eq:accuracy_ass}
\|z_m-H^{-1}v_m\|_H \leq \varepsilon \cdot \|H^{-1}v_m\|_H, 
\end{equation}
then
\begin{equation} \label{eq:norm_Vm_bound}
\|V_{m+1}\|_{H^{-1}} \leq \left( \frac{m+1}{1 - \varepsilon}\right)^{1/2}.
\end{equation}
\end{proposition}
\begin{proof}
Denote by $V_m \in \mathC^{n \times m}$ the matrix which contains the flexible Lanczos vectors $v_k$ as its $k$-th column. We have $x_m = x_0 + V_m\zeta_m$
with $\zeta_m$ satisfying the minimality condition \eqref{eq:fmr}. Thus,
using the flexible Lanczos relation \eqref{eq:Arnoldi2} and \eqref{eq:flex_initial_vector} we obtain
\begin{eqnarray*}
\|b-Ax_m\|_{H^{-1}} &=& \| V_{m+1}\beta_0 e_1 - V_{m+1} T_{m+1,m} \zeta_m\|_{H^{-1}} \\
&=&  \| H^{-1/2}V_{m+1}(\beta_0 e_1 -  T_{m+1,m} \zeta_m) \|_2 \\
&\leq & \| H^{-1/2}V_{m+1}\|_2 \cdot \|\beta_0 e_1 - T_{m+1,m} \zeta_m\|_2 \\
&=& \| H^{-1/2}V_{m+1}\|_2 \cdot \varrho_m,
\end{eqnarray*}
which is \eqref{eq:FMR_bound}. To obtain \eqref{eq:norm_Vm_bound} we now show that the $\ell_2$-norm of each column of $H^{-1/2}V_{m+1}$ is bounded by $\sqrt{1/(1-\varepsilon)}$. The inequality \eqref{eq:norm_Vm_bound} then follows as the $\ell_2$-norm of a matrix never exceeds its Frobenius norm. 

With $\nu_k := \|v_{k}\|_{H^{-1}} = \|H^{-1/2}v_k\|_2$ we have
\begin{eqnarray*}
\nu_k^2 = \langle v_k,z_k\rangle + \langle v_k, H^{-1}v_k-z_k \rangle 
 =  1 + \langle H^{-1/2}v_k,H^{-1/2}v_k-H^{1/2}z_k\rangle,
\end{eqnarray*}
from which we get, via the Cauchy-Schwarz inequality,
\[
| \nu_k^2-1 | \leq  \|H^{-1/2}v_k\| \cdot \|H^{-1/2}v_k-H^{1/2}z\| = 
\|v_k\|_H \cdot \|H^{-1}v_k - z_k\|_H.
\]
With assumption \eqref{eq:accuracy_ass} we therefore have
\[
| \nu_k^2-1 | \leq \varepsilon \nu_k^2,
\]
which is equivalent to
\[
\frac{1}{1+\varepsilon} \leq \nu_k^2 \leq \frac{1}{1-\varepsilon}.
\]
\end{proof}

There are two important observations regarding Proposition~\ref{prop:bound}.
First, practically, we cannot monitor $H$-norms of the error, so inequality \eqref{eq:accuracy_ass} cannot be checked directly. However, a relative reduction of the ($\ell_2$-norm of the) residual, which can be 
measured, typically results in a similar reduction of the error. Moreover, in 
the case that we use the CG method to approximate $H^{-1}v_k$, it is precisely
the $H$-norm of the error which is made as small as possible by the CG iterates. We refer to \cite{Sarkis.Szyld.07} for 
a further discussion of this and related aspects. 

Secondly, the quantity $\varrho_m$ is available in the algorithm at no additional cost. The situation is similar as in GMRES \cite{Saad03} or QMR \cite{FreNac91}: To solve the least sqares problem \eqref{eq:fmr}, we update the $QR$-factorization $T_{m+1,m} = Q_{m+1,m+1}R_{m+1,m}$ in each step. The value $\varrho_m$ is the last, i.e.\ the $(m+1)$st, entry of $Q_{m+1,m}^*\|r_0\|_H e_1$, and this value is updated to $\varrho_{m+1}$ using the Givens rotation which updates $Q_{m+1,m}$ to $Q_{m+2,m+1}$.  

We further note that the proof of Proposition~\ref{prop:bound} also shows that in the case where we solve for $H$ exactly, we find $\varrho_m = \| b-Ax_m \|_{H^{-1}} $ as $\|V_m\|_{H^{-1}} = 1$.

\section{Numerical results} \label{numerics:sec}
In this section we report results for several examples which have been considered in the literature. The general setup for all experiments is as follows: We report residual norm plots for both, the FMR and the FGAL method. In each iteration,  the systems with the hermitian positive definite part $H$ are solved approximately using the CG method, asking for the initial $\ell_2$-norm of the residual to be reduced by a factor of $\varepsilon_{\cg }$. We typically vary $\varepsilon_{\cg }$ from $10^{-1}$ up to $10^{-12}$. In light of \cref{prop:bound} and the discussion preceding it,  the FGAL and FMR iteration is stopped once the $H^{-1}$-norm of the initial residual is reduced by $\varepsilon_{f} = 10^{-12}$. The case $\varepsilon_{\cg }=10^{-12}$ in the approximate solves may thus actually be considered as a equivalent to exact solves in our context. 

The residual norm plots report (relative) approximate $H^{-1}$-norms of the residual, i.e.\ the quantity $\sqrt{\widehat{r}^*r}$ where $r$ is the current residual and $\widehat{r}$ is the result of CG when solving $H\widehat{r} = r$ with initial guess $0$ and a required reduction of the residual $\ell_2$-norm of $\varepsilon_{\cg }$. This number can be somewhat off the exact $H^{-1}$-norm;
the advantage is that it is computed with an effort comparable with what we need to update an iterate. When the iterations stopped, we double-checked that the targeted decrease of the $H^{-1}$-norm was indeed achieved by
recomputing $\sqrt{\widehat{r}^*r}$ for the final residual with now $\varepsilon_{\cg } = 10^{-12}$ in the CG method which computes $\hat{r}$. Finally, as an illustration, our plots also give the bound on the $H^{-1}$ norm of the residuals of the FMR iterates according to \cref{prop:bound}. Let us stress that, as opposed to the approximate $H^{-1}$-norms, these bounds can be computed at virtually no extra cost, so that in a production setting we could base a stopping criterion for the FMR iterations merely on these bounds.

%Our examples focus on implicit time integration for dissipative Hamiltonian ODEs, but

We start reconsidering the example in \cite{Widlund78}.

\paragraph{Example 1} Consider the convection-diffusion equation
\begin{eqnarray*}
-\Delta u(x,y) + au_x(x,y) = b(x,y),
\end{eqnarray*}
on $\Omega = [0,1] \times [0,1]$ with Dirichlet boundary conditions. 
We discretize on a uniform square mesh with $127 \times 127$ interior grid points, using central finite differences in both, $\Delta u$ and $u_x$,  and choose $a = 10^4$. The right hand side $b$ was chosen to have random entries.

\begin{figure} 
\includegraphics[width=0.49\textwidth]{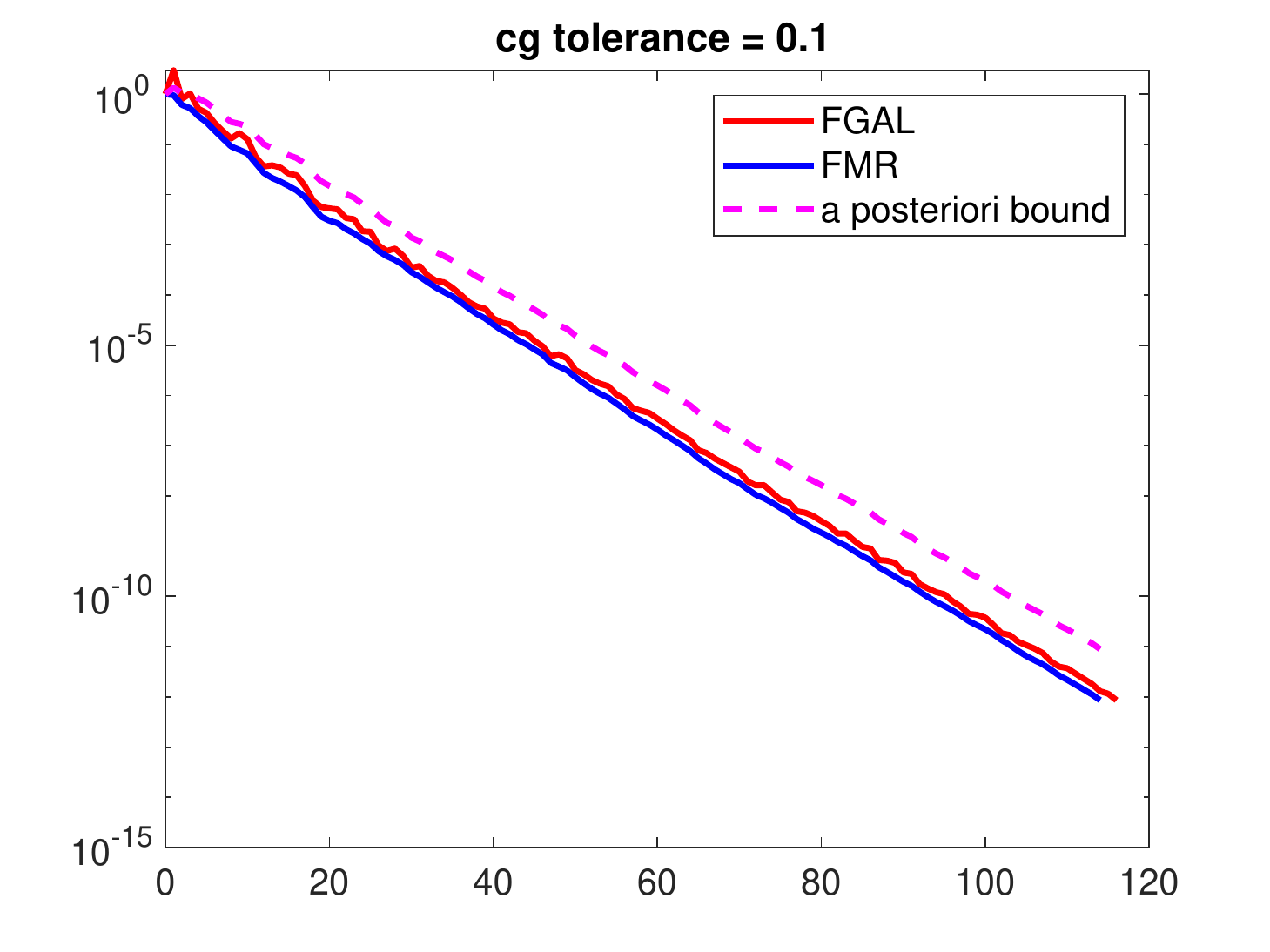}
\includegraphics[width=0.49\textwidth]{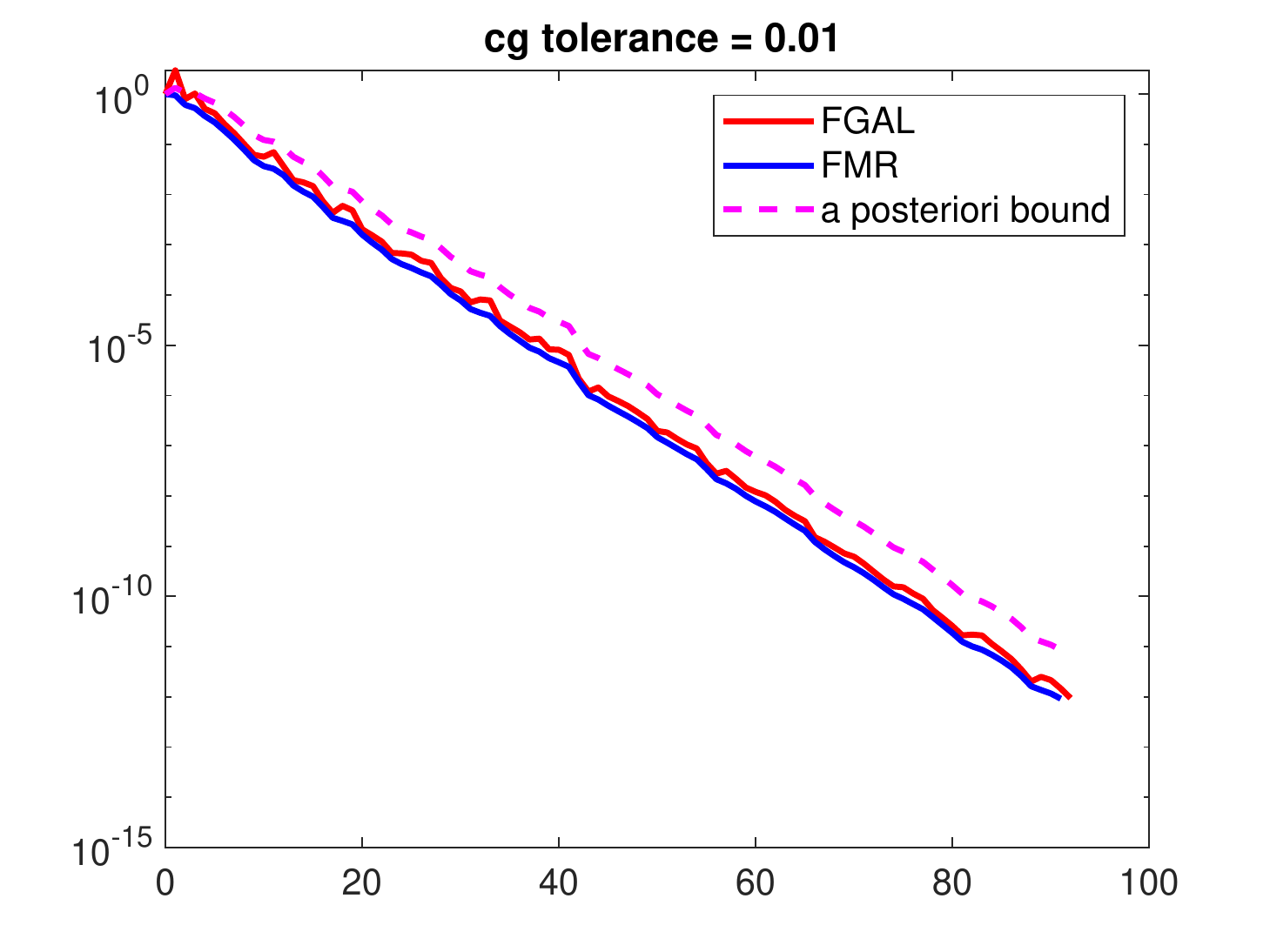}
\includegraphics[width=0.49\textwidth]{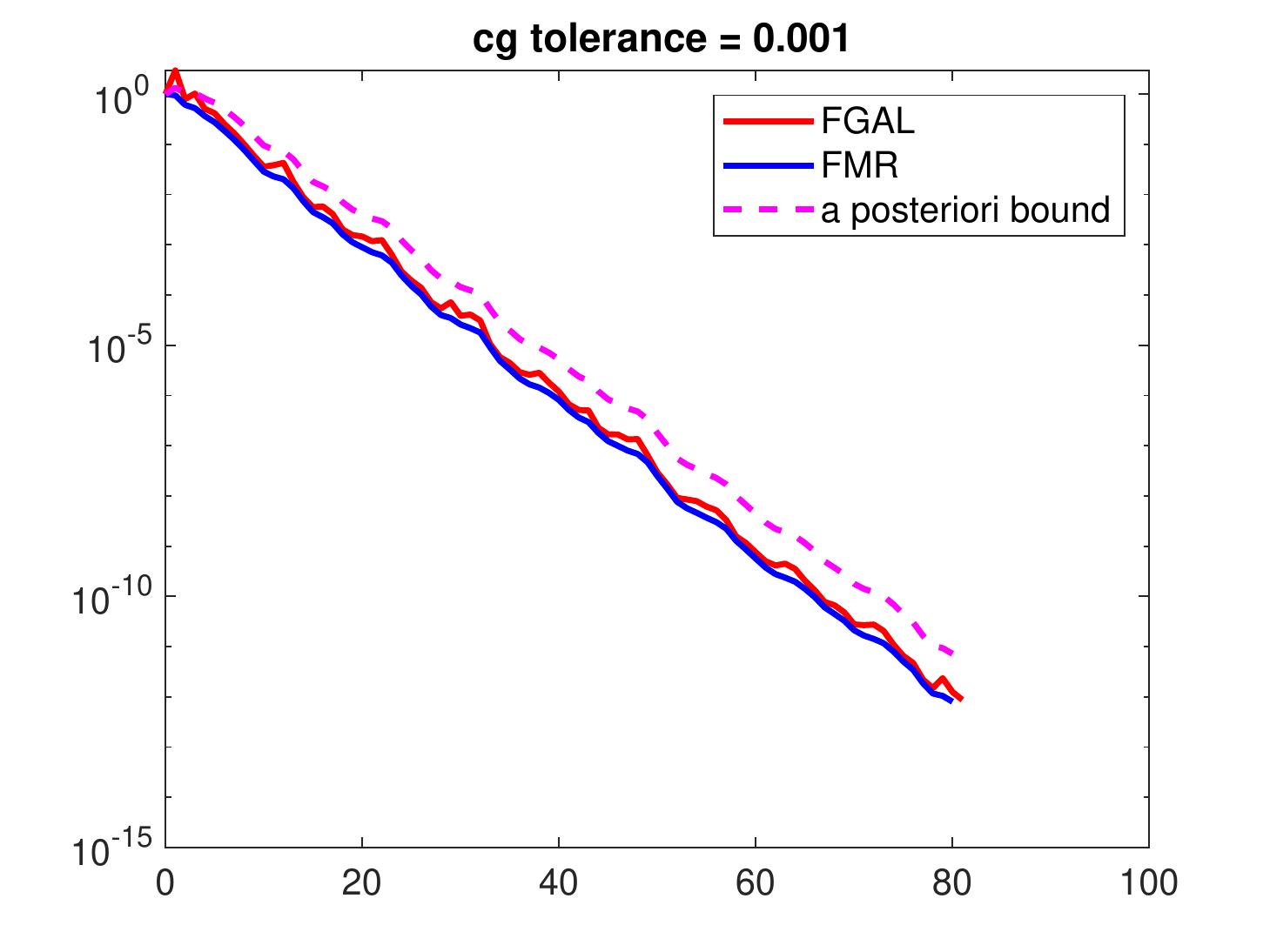}
\includegraphics[width=0.49\textwidth]{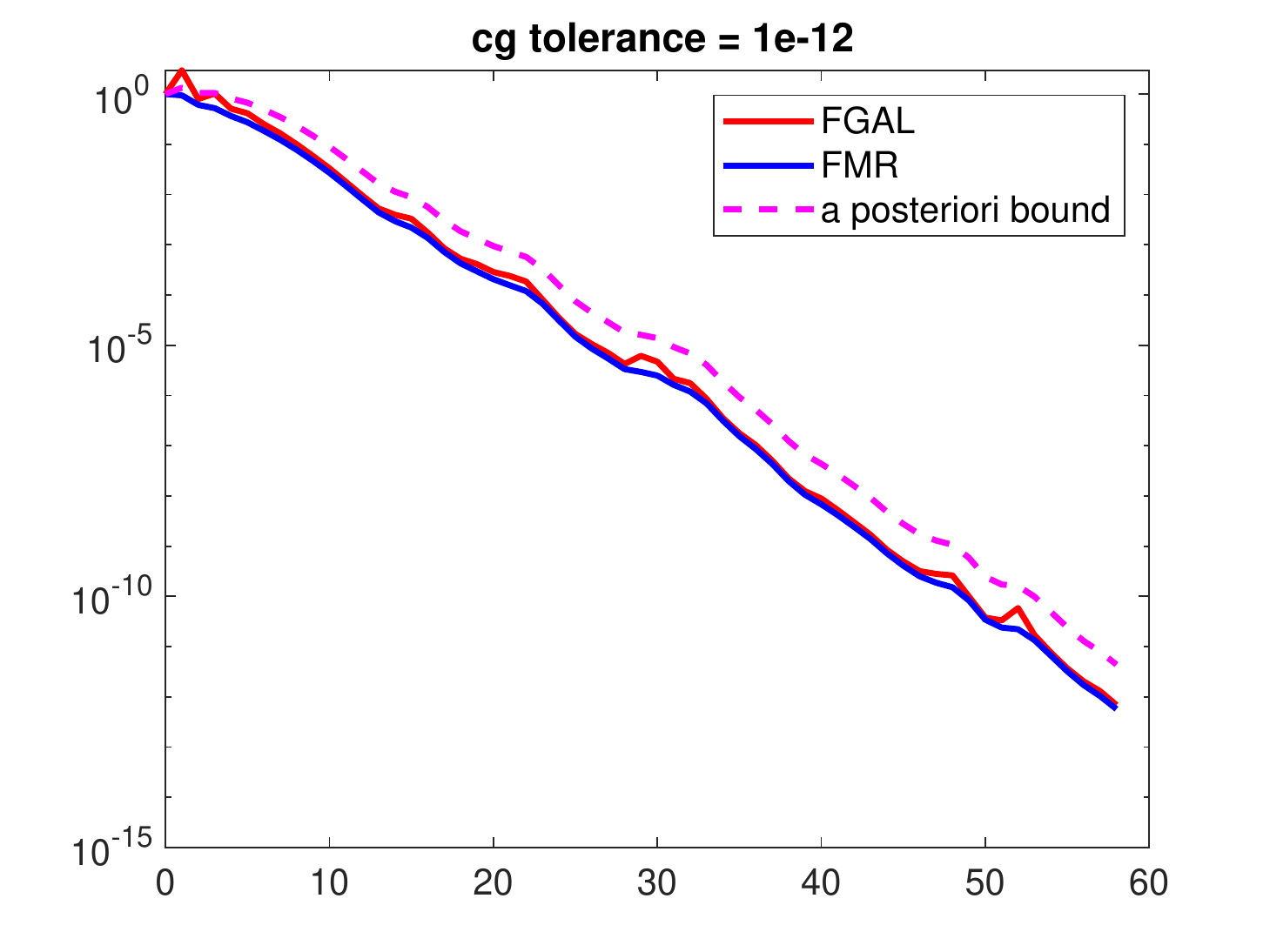}
\caption{Convergence plots for the convection-diffusion equation, Example~1, four choices for $\varepsilon_{\cg }$ \label{fig:widlund}} 
\end{figure}

Figure~\ref{fig:widlund} shows convergence plots for four different choices of $\varepsilon_{\cg }$. We see that for all $\varepsilon_{\cg }$ the
residual norms of the FGAL and the FMR iterates differ only marginally, with the FMR residuals, as is to be expected due to their (approximate) minimization property, being (slightly) smaller and decreasing in a smoother manner. We also see that the bound from \cref{prop:bound} captures the actual convergence behavior quite accurately. With $\varepsilon_{\cg } = 10^{-1}$, we already need only about twice as many iterations as with ``exact'' solves, $\varepsilon_{\cg } = 10^{-12}$, while an approximate solve with $\varepsilon_{\cg } = 10^{-1}$ requires about 50 CG iterations on average as opposed to about 470 CG iterations on average for $\varepsilon_{\cg } = 10^{-12}$.

\paragraph{Example~2} This is the matrix treated in section 6.2 of \cite{GueLieMehSzy21}. It is based on a discretized Stokes problem obtained by using the $Q_1$-$Q_1$ finite element discretization of the unsteady channel domain problem in IFISS \cite{IFISS}. This example was treated for two different grid parameters in \cite{GueLieMehSzy21}. We only consider the larger of the two, i.e\ grid parameter 10, resulting in a linear system of size $n =  3,151,875$, and we took a step size of $\tau/2 = 10^{-3}$ in the midpoint Euler rule. 

\begin{figure} 
\centerline{
\includegraphics[width=0.5\textwidth]{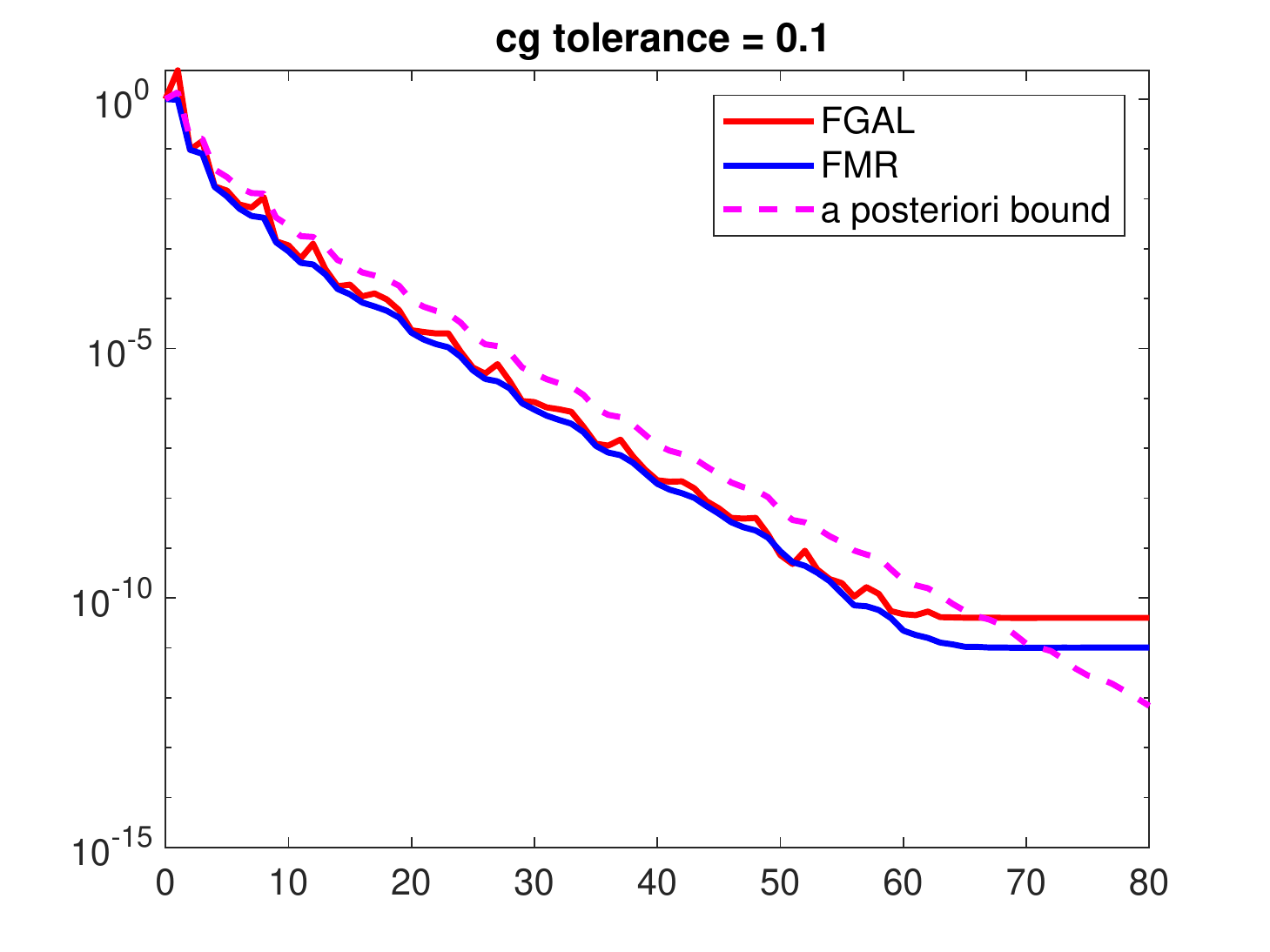}
\includegraphics[width=0.5\textwidth]{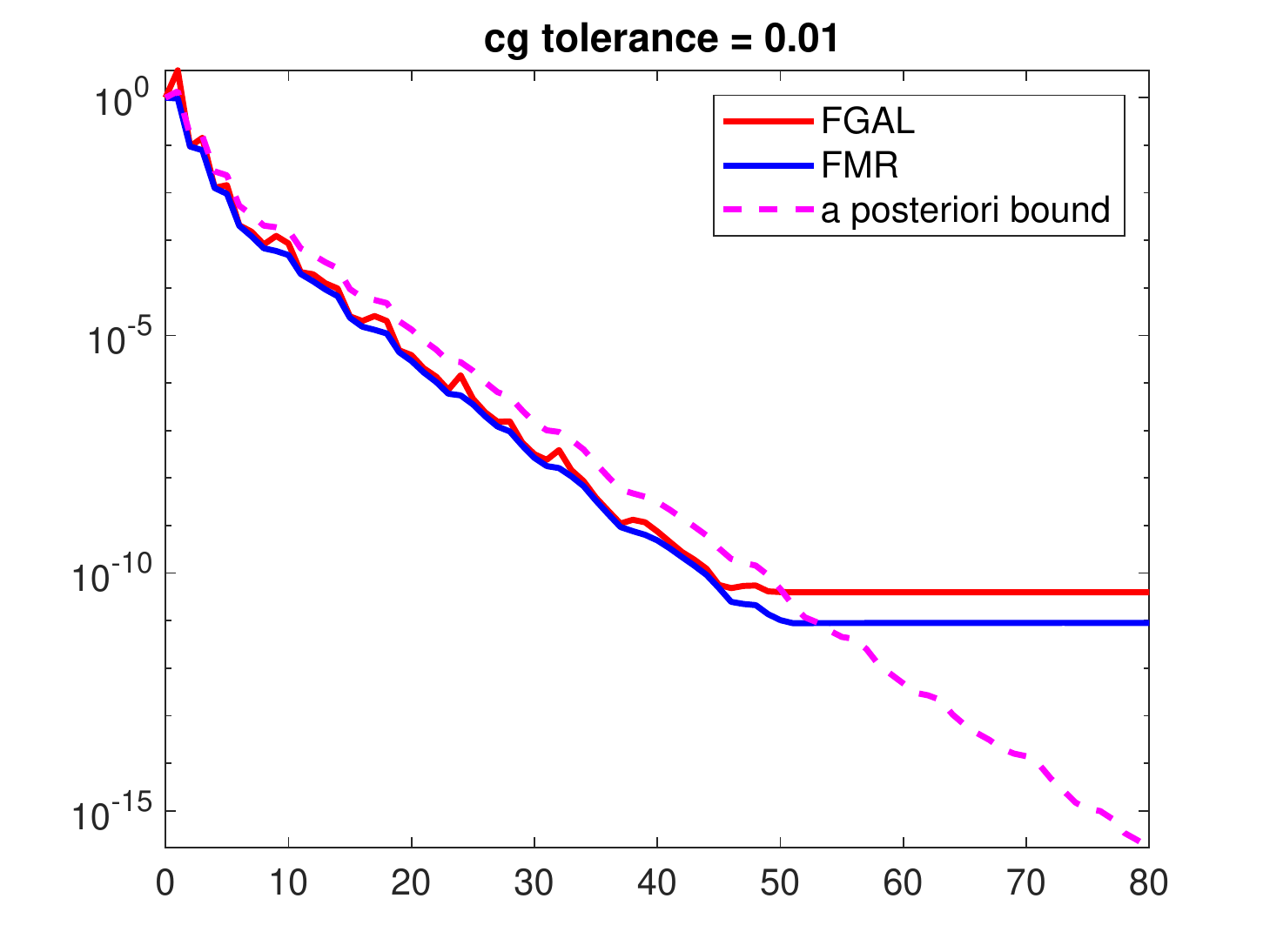}
}
\centerline{
\includegraphics[width=0.5\textwidth]{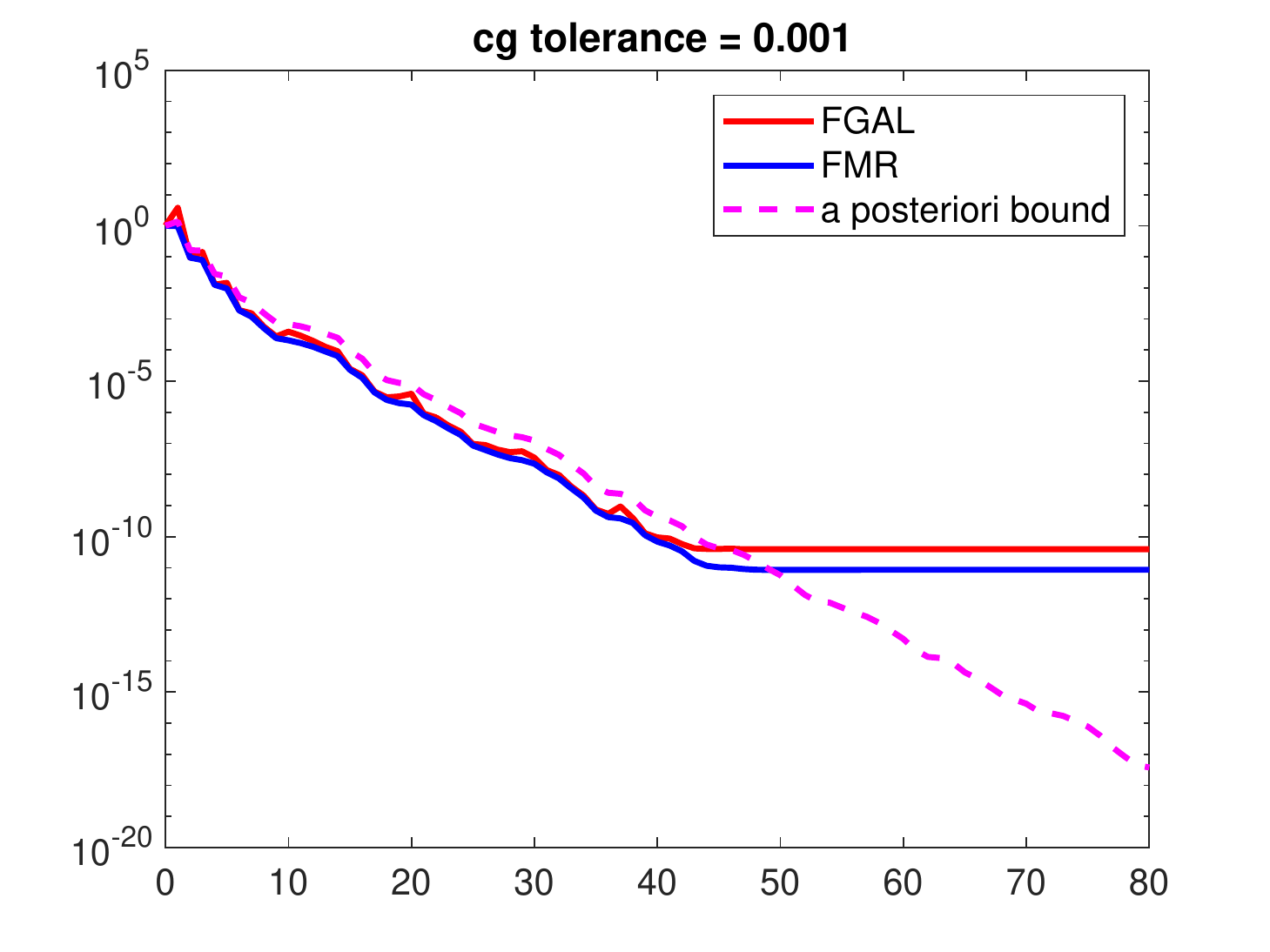}
\includegraphics[width=0.5\textwidth]{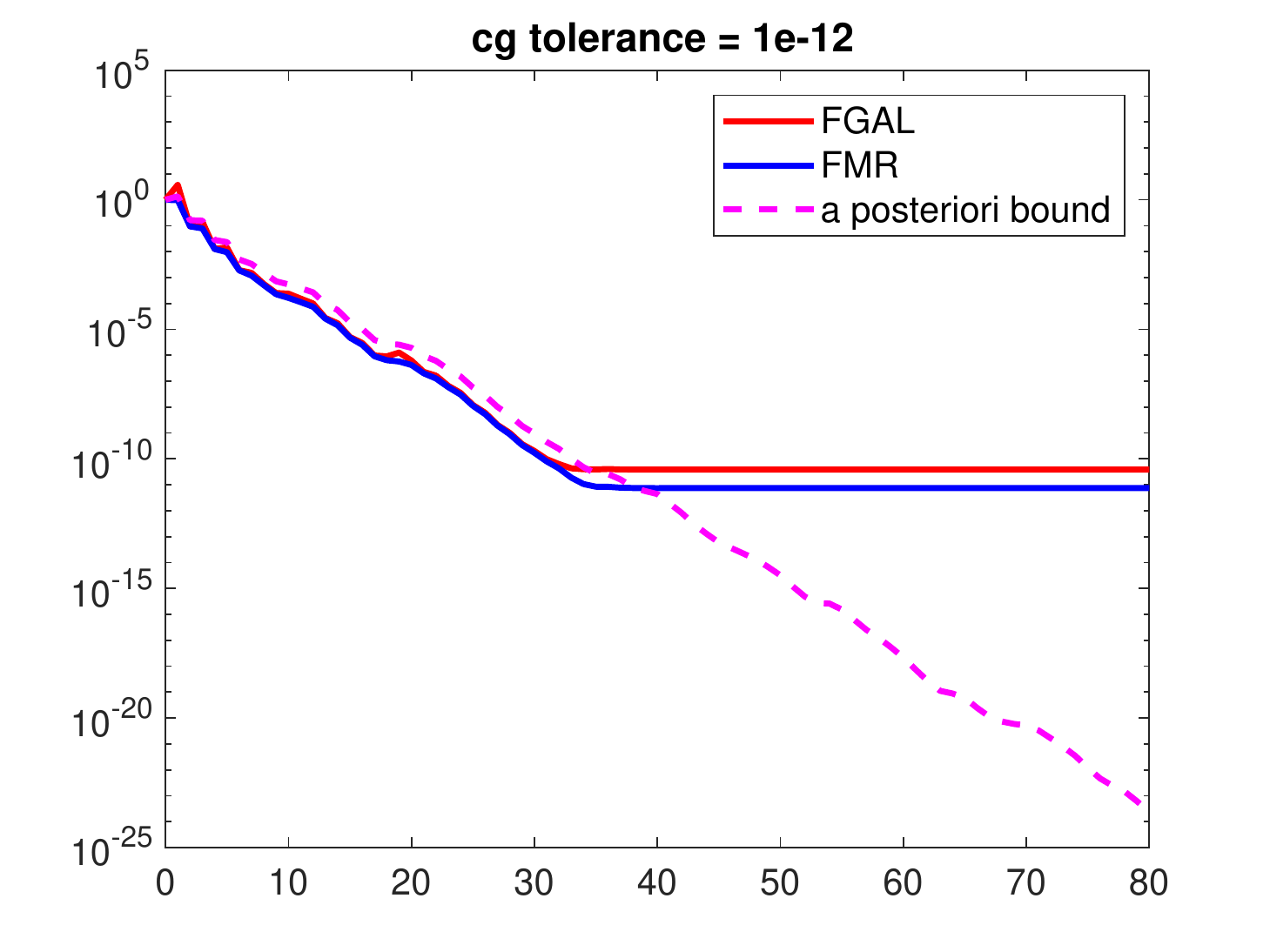}
}
\caption{Results for the unsteady channel domain Stokes problem, Example 2, $\tau/2 = 10^{-3}$, four choices for $\varepsilon_{\cg }$ \label{fig:stokes_verybig}}
\end{figure}

The results are reported in Figure~\ref{fig:stokes_verybig}. To solve the systems with $H$ we use preconditioned CG with the ILU(0) preconditioner. In the solves for $H$, every about 12 preconditioned CG iterations reduce the norm of the residual by another factor of $10^{-1}$. 

We observe that the dependence on the CG tolerance $\varepsilon_{\cg}$ is similar as in the previous example, with $\varepsilon_{\cg} = 10^{-1}$ already 
requiring just about half times more the number of iterations than with ``exact'' solves ($\varepsilon_{\cg } = 10^{-12}$). FGAL and FMR behave similarly with FGAL showing a slightly more oscillatory convergence behavior than before. This example also shows stagnation once the residual norm is reduced by a factor of approximately $10^{-10}$ for FGAL and $10^{-11}$ for FMR. Note that here as in all the other plots we always report the norms of ``computed'' residuals, i.e.\ residuals which were explicitly computed from the current iterate in each iteration. We attribute the observed stagnation to the relatively high condition number $\kappa$ of the matrix which makes us expect that the reported residual norms, which rely on explicitly computed residuals, will not go beyond roughly $\kappa \cdot \varepsilon_{\text{mach}}$, with $\varepsilon_{\text{mach}} \approx 10^{-16}$ being the machine roundoff unit. The reported residual \textit{bounds} are based on updated quantities, which is why they continue to decrease beyond the point where the computed residuals stagnate.

Interestingly, the numerical computations in \cite{GueLieMehSzy21} used only an \textit{incomplete} Cholesky factorization to solve the systems with $H$---the incomplete Cholesky drop tolerance in Matlab's \texttt{ichol} was set to the relatively small value of $10^{-9}$. Thus, the experiments reported in \cite{GueLieMehSzy21} perform approximate solves only, while the algorithms used there are not flexible. This might be the reason why the experiments in \cite{GueLieMehSzy21}, as the authors say, could not use larger values than $10^{-3}$ for $\tau/2$. 
\cut{In Figure we report results for the flexible methods for $\tau/2 = 10^{-3}$ and $\tau/2 = 10^{-2}$, again using different values for $epsilon_{\cg }$.
}

The remaining two examples are linear systems that arise when using the implicit midpoint Euler rule when integrating a DOE resulting from a port Hamiltonian modelling. These DOEs have the form
\begin{equation} \label{DAE:eq}
E\dot{u}(t) = (J-R)u(t) + f(t), \enspace u(0) = u_0,
\end{equation}
where $E$ is hermitian and positive definite, $J$ is anti-selfadjoint and $R$ is hermitian and positive semidefinite. With a step size of $\tau > 0$, the linear systems to solve at each time step are then of the form
\[
A x = b \enspace \mbox{ with } A = E + \frac{\tau}{2}(R-J),
\]
and we have $H = E + \frac{\tau}{2}R$, positive definite, and $S = -\frac{\tau}{2}J$. Note that implicit integration methods are mandatory in this situation if one wants to preserve the dissipation of energy; see \cite{Mehrmann_2019}.

\begin{figure} 
\centerline{
\includegraphics[width=0.49\textwidth]{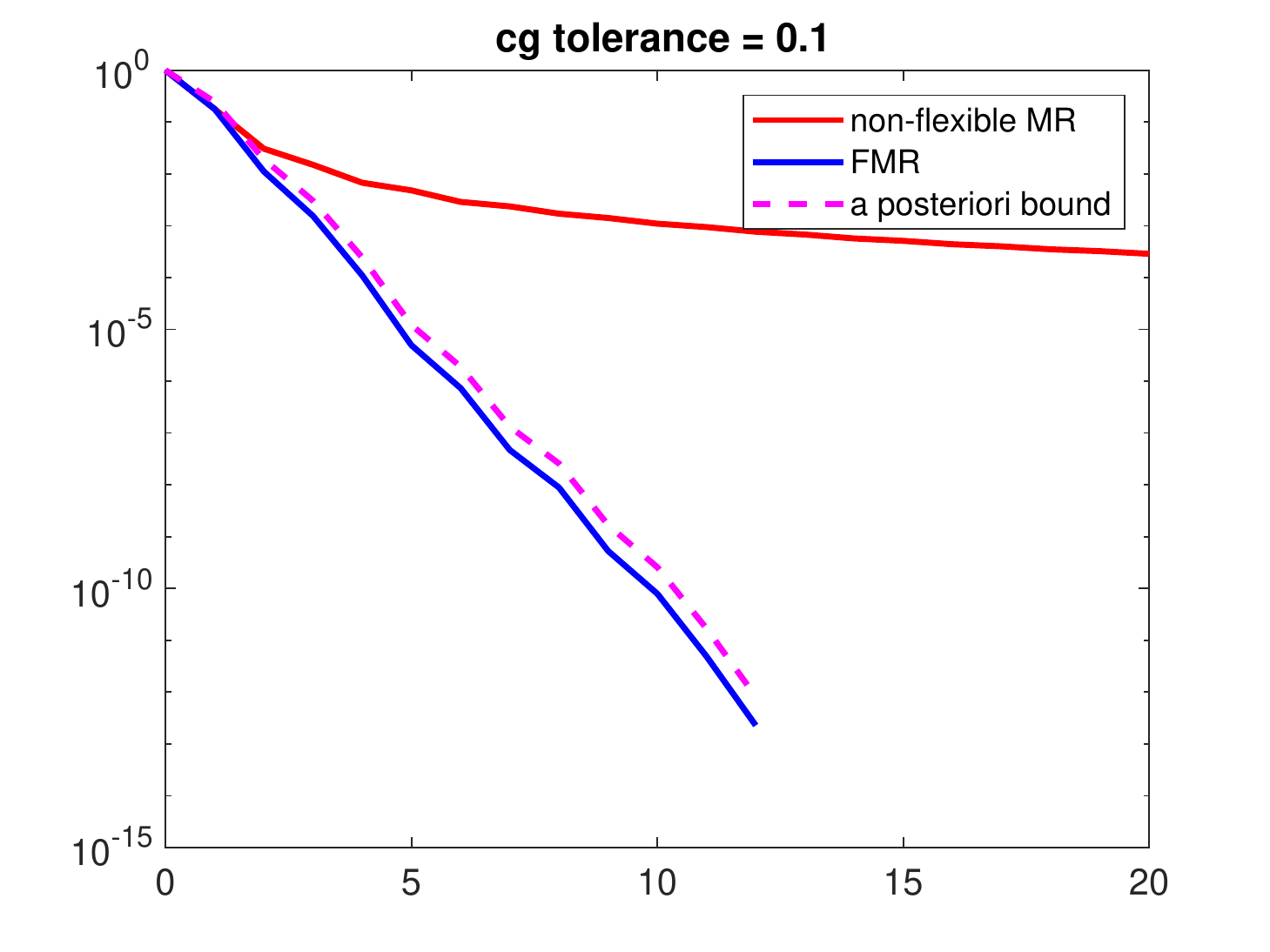}
\includegraphics[width=0.49\textwidth]{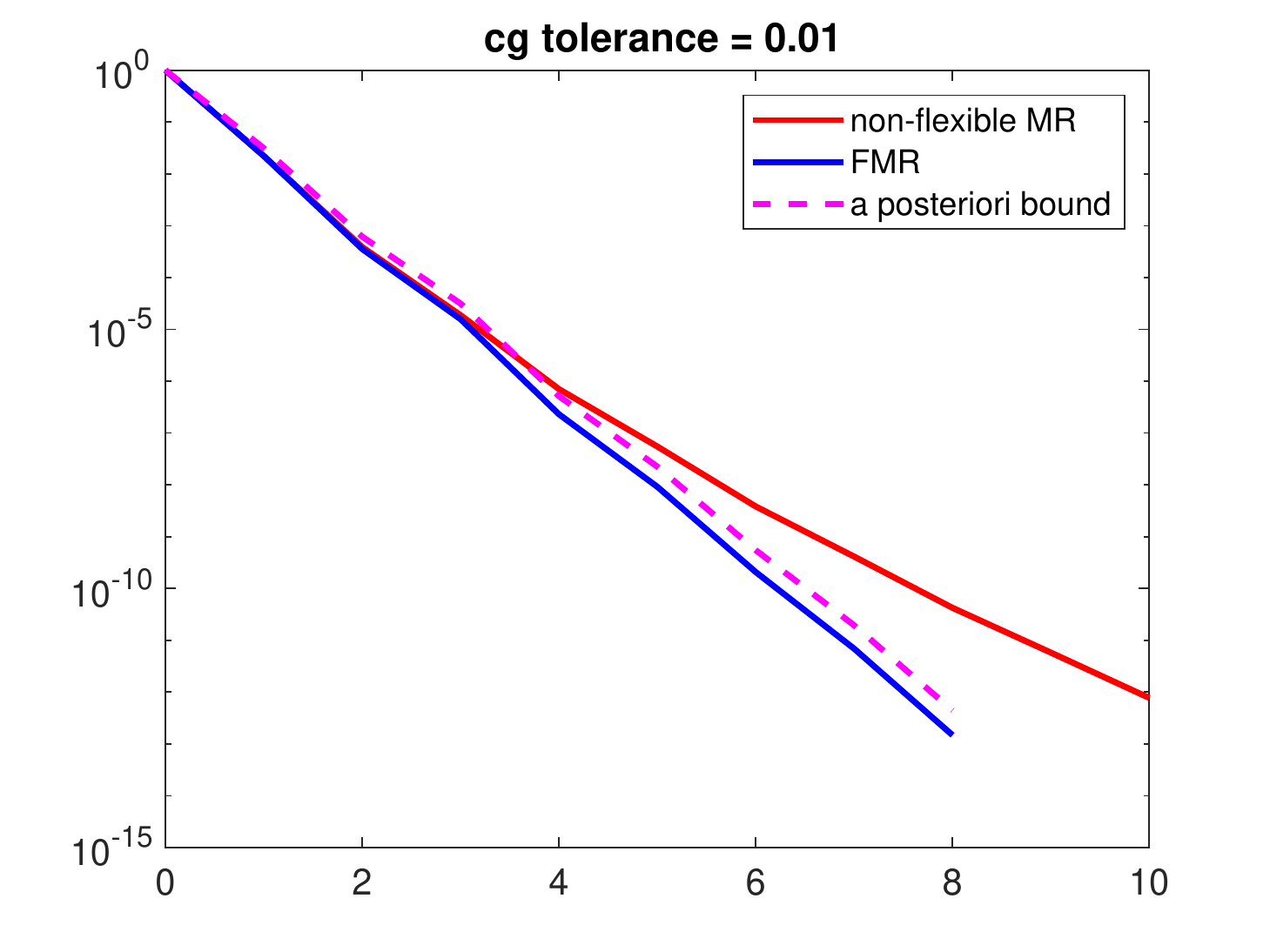}
}
\centerline{
\includegraphics[width=0.49\textwidth]{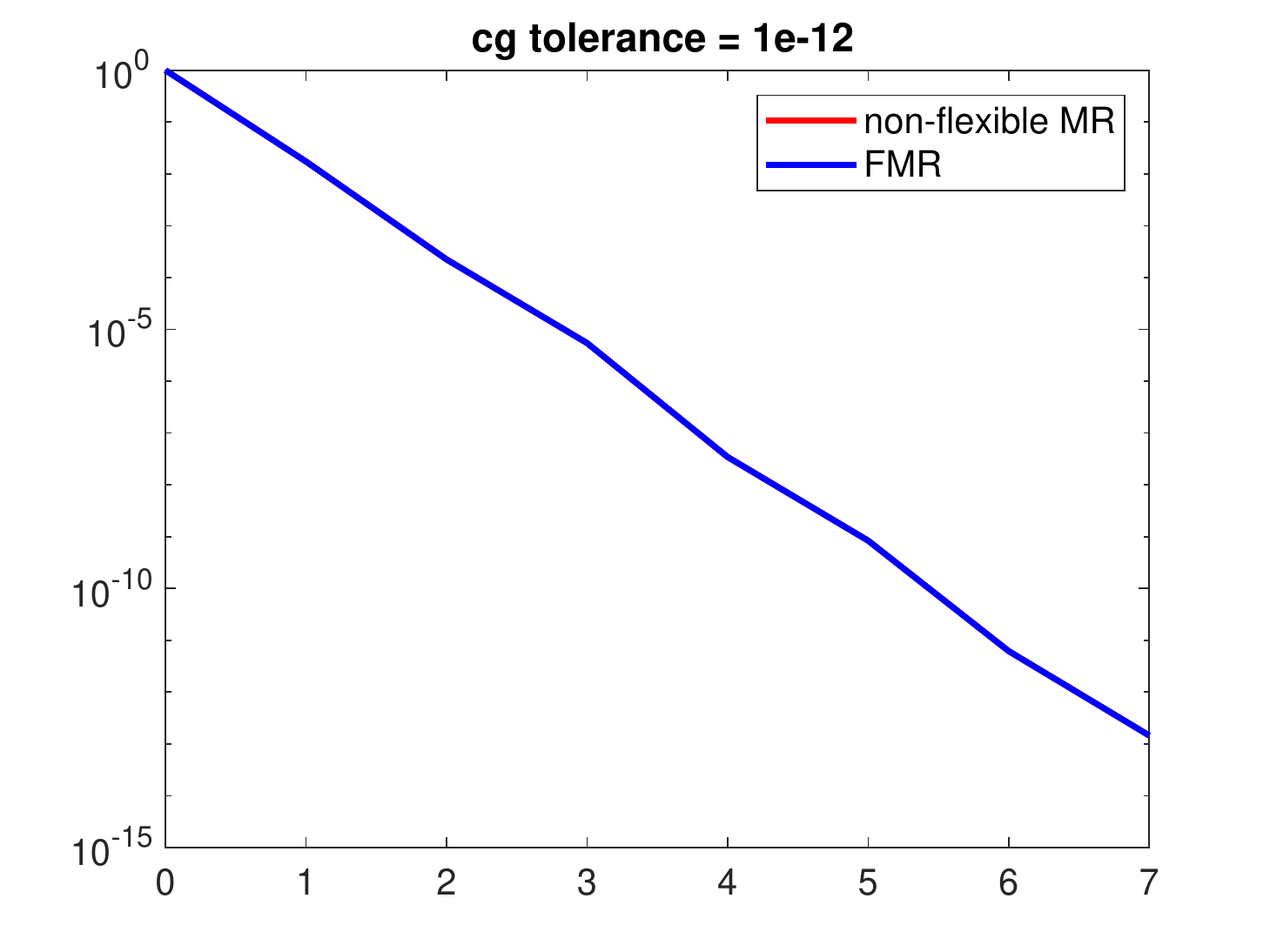}
\includegraphics[width=0.49\textwidth]{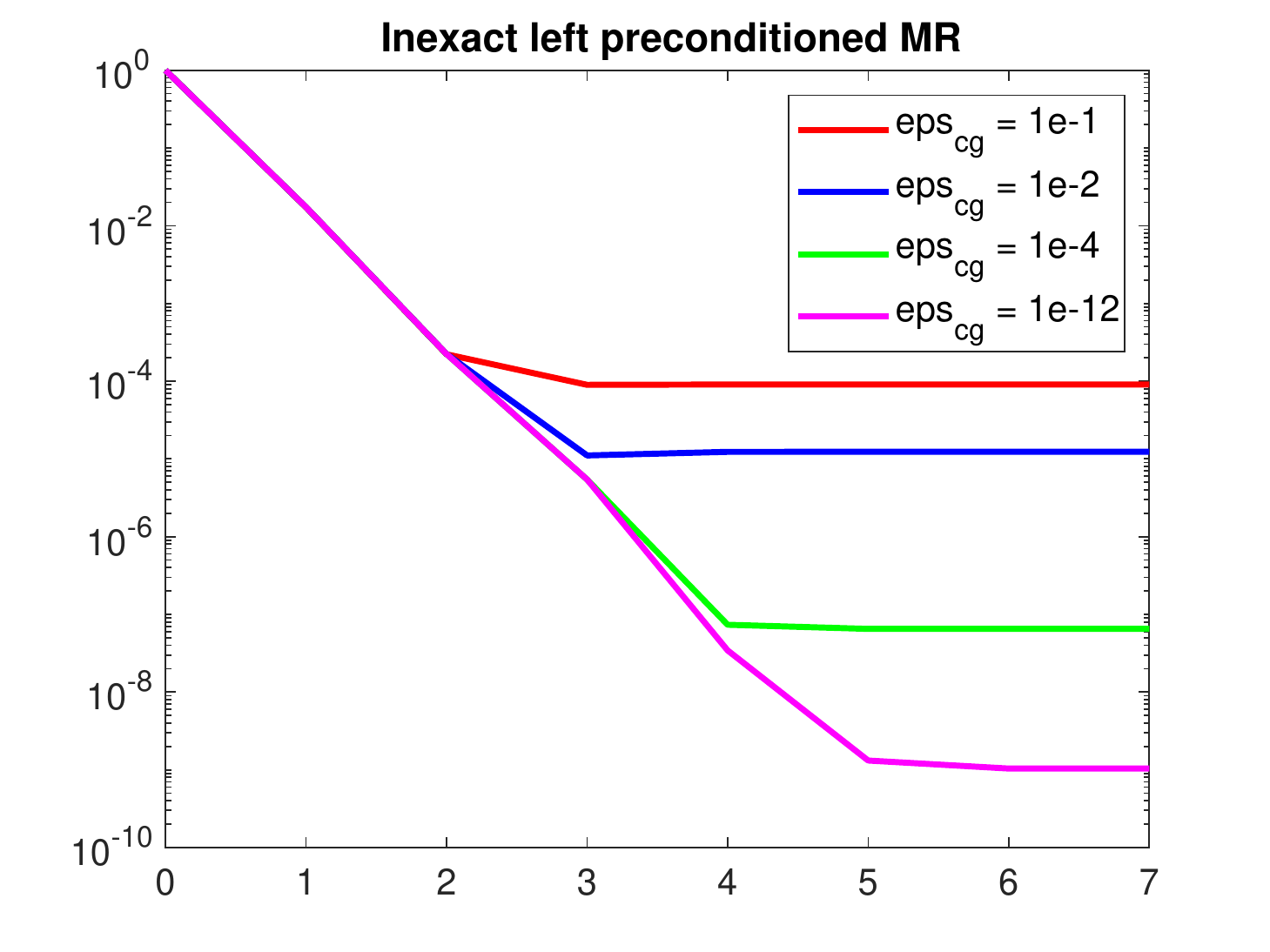}
}
\caption{Convergence plots for the spring mass system of dimension $10^6$, Example~3, $\tau/2 = 10^{-1}$, three choices for $\varepsilon_{\cg }$. Bottom right plot is for non-flexible $\ell$MR with inexact solves for $H$. \label{fig:spring1Mio}} 
\end{figure}

\paragraph{Example~3} We take the homonomically constrained spring mass system from \cite{MehSty05}, which is part of the port Hamiltonian system benchmark collection\footnote{available at \url{https://algopaul.github.io/PortHamiltonianBenchmarkSystems/}}.
%\cite{pHS_bench}. 
We refer to \cite{GPBS2012} and \cite{GueLieMehSzy21} on how the matrices $E$, $R$ and $J$ arise from the descriptor system formulation of the system. We took the system with $n=2\cdot 10^6$ as it was also considered in \cite{GueLieMehSzy21} and $\tau/2 = 0.1$.  The numerical results are given in Figure~\ref{fig:spring1Mio}. For all choices of $\varepsilon_{\cg }$, the norms of the residuals of FGAL and FMR are so close that they are indistinguishable in the plots, which is why we report only the results for FMR. The bounds from \cref{prop:bound} again match the convergence behavior quite exactly. Similar to the previous two examples, we see that with $\varepsilon_{\cg } = 10^{-1}$ we need less than twice as many iterations than with ``exact'' solves ($\varepsilon_{\cg } = 10^{-12}$), and this number of iterations is this time exceeded by just 1 for $\varepsilon_{\cg } = 10^{-2}$. We need between 4 and 6 (unpreconditioned) CG iterations with $H$ to reduce the initial residual by the factor $\varepsilon_{\cg} = 10^{-1}$ in each step of $\ell$MR, and between 43 and 49 to reduce it by the factor $\varepsilon_{\cg} = 10^{-12}$. 

In the three plots for $\varepsilon_{\cg } = 10^{-1}, 10^{-2}$ and $10^{-12}$ we in addition plot the $H$-norms of the residuals in a variant of FMR which differs from FMR only by the fact that in the underlying preconditioned Lanczos algorithm, Algorithm~\ref{flexlanczos:alg}, we do not compute $\gamma_k$ but rather put $\gamma_k = \beta_{k-1}$. This method is termed ``non-flexible MR''. The plots show that for the larger values of $\varepsilon_{\cg }$ this modification results in a significantly to dramatically delayed convergence. 

In a similar spirit, the bottom right plot of Figure~\ref{fig:spring1Mio}
show the $H^{-1}$-norm of the residual for the non-flexible $\ell$MR method of 
Rapoport (left preconditioning), where we solve for $H$ with CG at different accuracies. Not really surprisingly, and convincingly, the plot illustrates that this ``inexact'' version of $\ell$MR stagnates sooner the less accurately we solve for $H$. A similar behavior can be observed for non-flexible inexact $\ell$GAL, for which we do not report results here.
In contrast, the flexible methods reach high accuracies for the overall system even when $\varepsilon_{\cg }$ is large.

%%%%%%%%%%%%%%%%%%%%%%%%%%%%%%%%%%%%%%%%%%%%%%%%%%%%%%
\paragraph{Example~4} We take the coupled electro-thermal DAE system from \cite{coupledET}, from which we treat the thermal part after applying the decoupling described in \cite{coupledET}. This again yields an ODE of the form \eqref{DAE:eq}. 
%resulting thermal subsystem has the form
%\begin{equation*}
%    E_T x' + A_T x = Bu(t).
%\end{equation*}
Similar to the previous example, we consider the linear systems arising when using the implicit midpoint Euler method for time integration. 
%resulting in the linear system $Ax=b$ with $A=E_T + \frac{\tau}{2} A_T$ and we have $H = (A+A^\top)/2$ is hermitian positive definite and $S=A-H$ is anti-selfadjoint. 
The system has dimension $n=7.2 \times 10^5$, and we take $\tau /2 = 10^{-2}$. Since the positive definite diagonal matrix $E$ is highly ill-conditioned, we scaled the system from left and right with the square root of the inverse of $E$. 

\begin{figure} 
\centerline{
\includegraphics[width=0.5\textwidth]{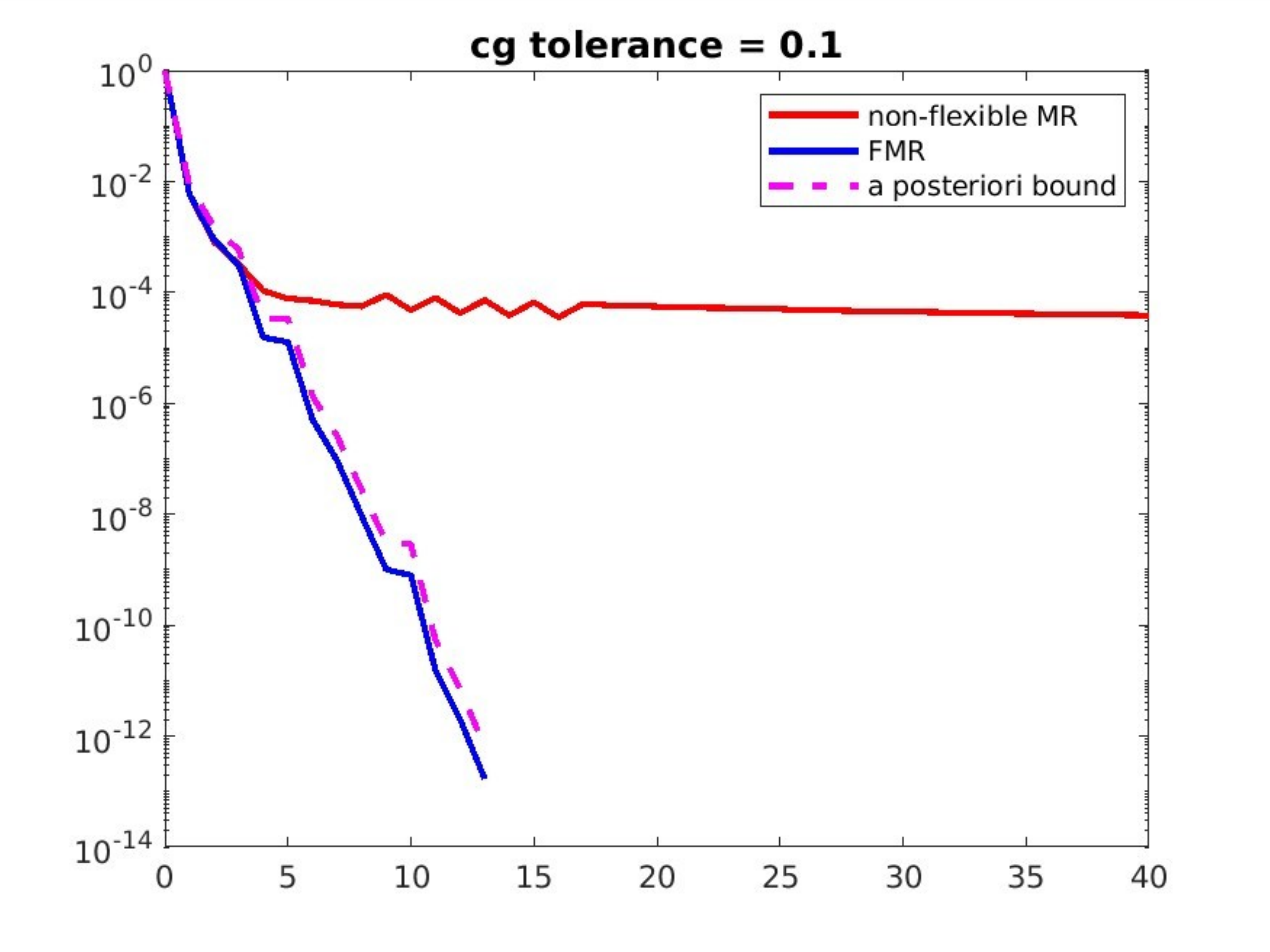}
\includegraphics[width=0.5\textwidth]{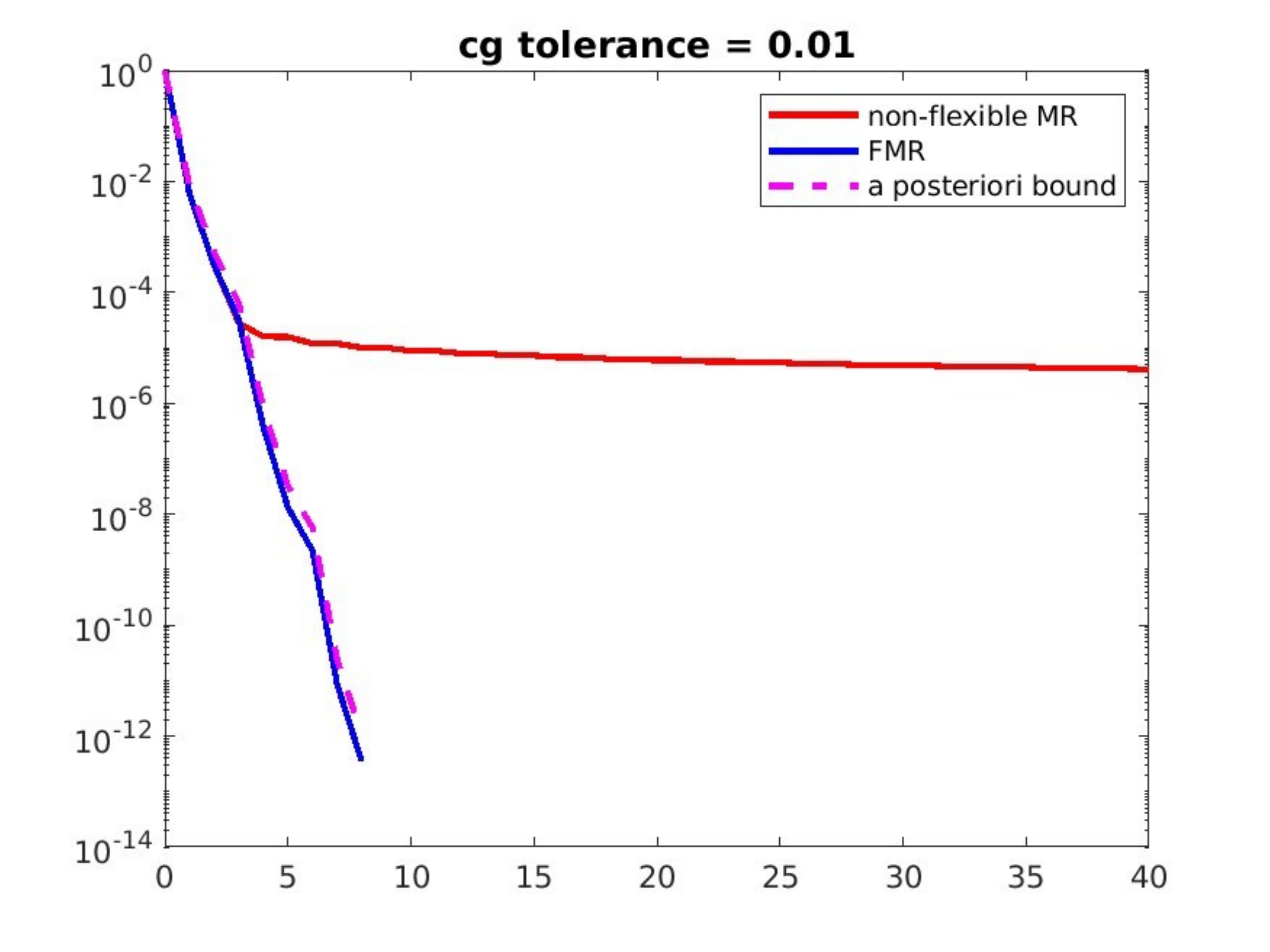}
}
\centerline{
\includegraphics[width=0.5\textwidth]{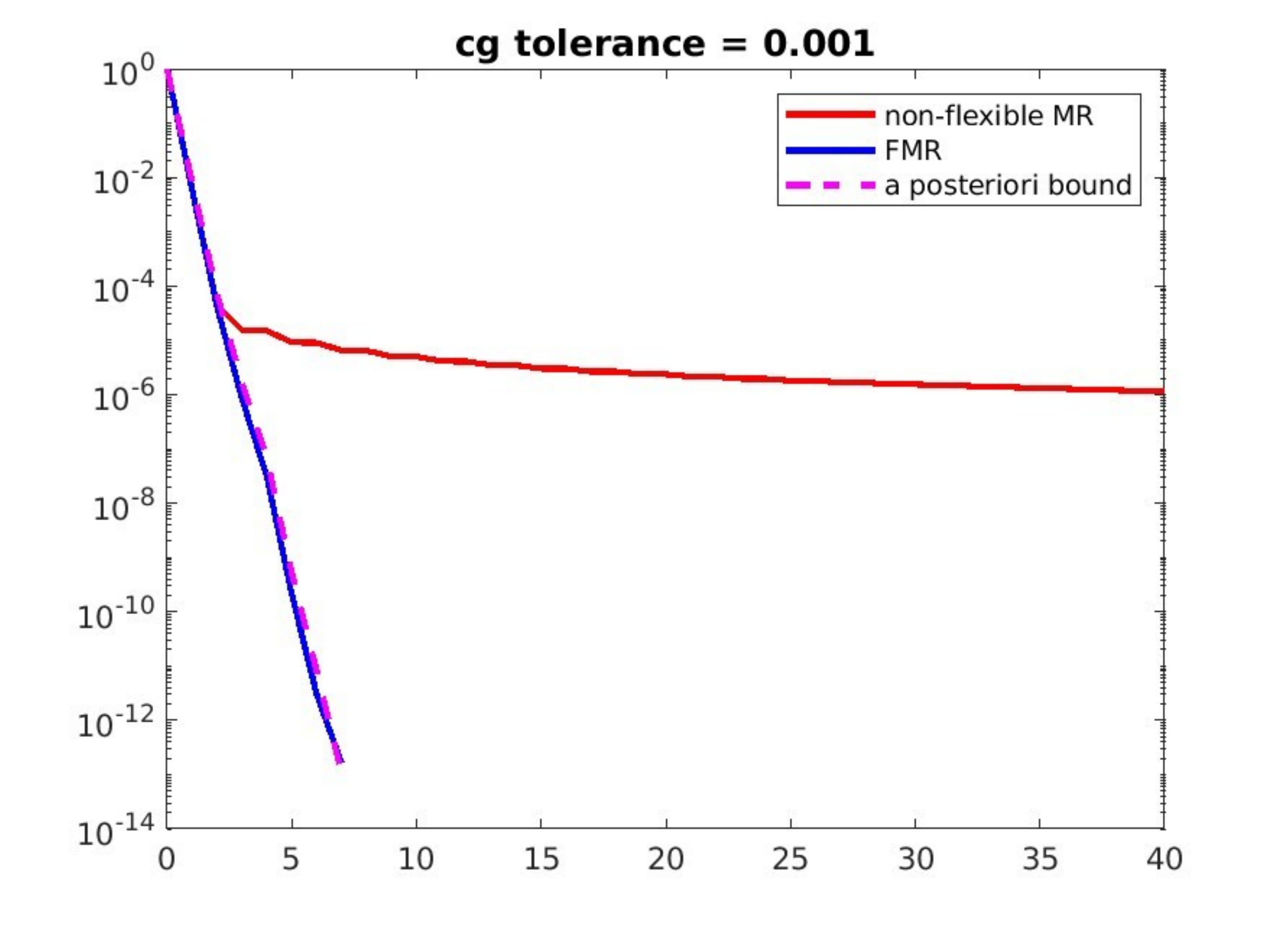}
\includegraphics[width=0.5\textwidth]{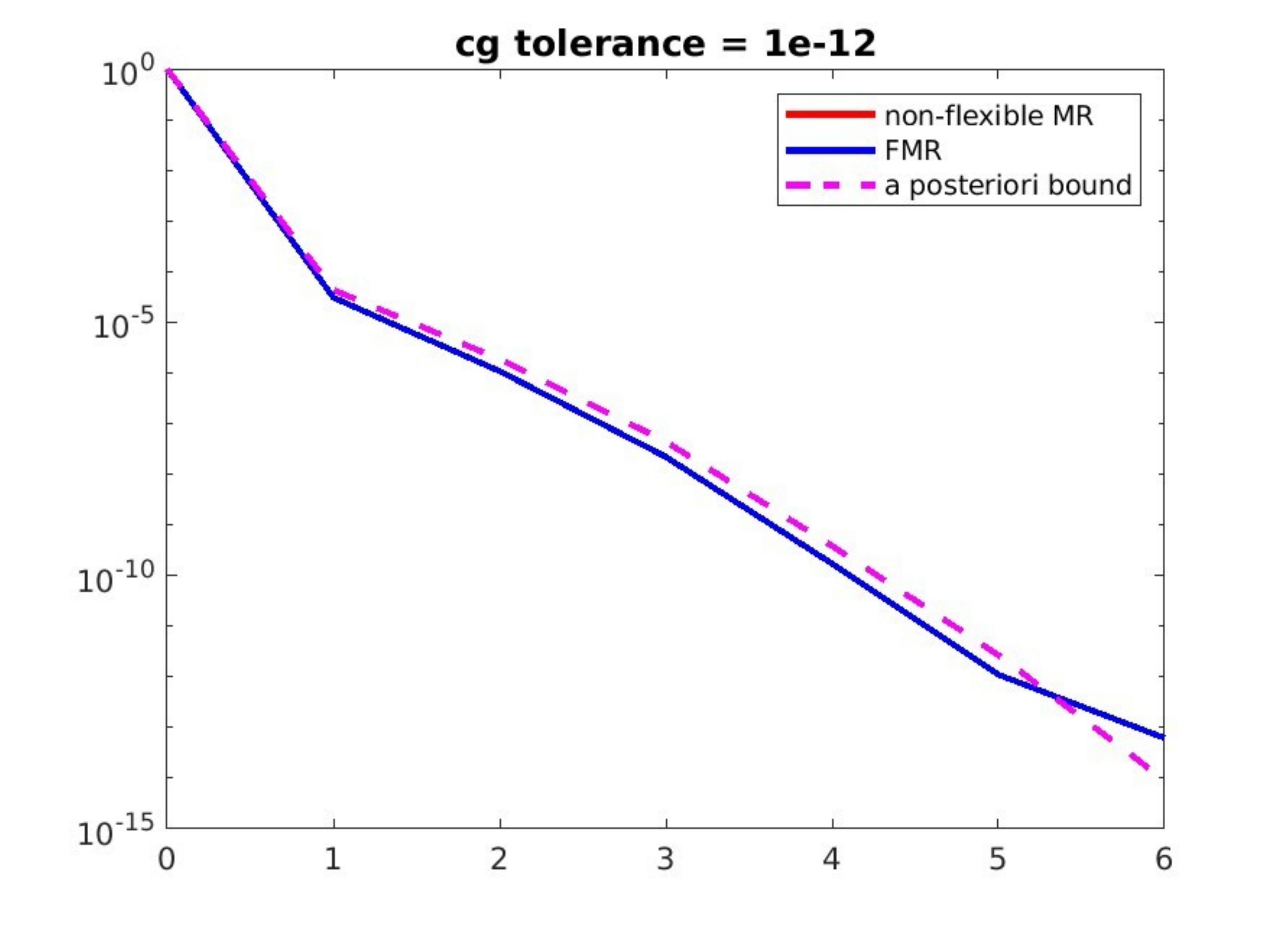}
}
\caption{Results for the coupled electro-thermal problem in Example 4 with $\tau/2 = 10^{-2}$ and four choices for $\varepsilon_{\cg }$ are considered. \label{fig:thermal}}
\end{figure}

The results are reported in Figure \ref{fig:thermal}. To solve the systems with $H$ we use preconditioned CG with a threshold dropping incomplete Cholesky preconditioner with drop tolerance $10^{-4}$. This results in a very efficient preconditioner, since on average we need just slightly less than two iterations to decrease the residual by one order of magnitude. Similarly to the previous example, FMR and FGAL perform very similarly and would be indistinguishable in the plots, which is why we only report the results for FMR. As before, we observe that FMR (and FGAL)  
require about twice the iterations with $\varepsilon_{\cg} = 10^{-1}$   than with ``exact'' solves ($\varepsilon_{\cg } = 10^{-12}$). Both flexible methods again achieve small final residuals even when large values of $\varepsilon_{\cg }$ are used.

For non-flexible MR, the iterations stagnate quite early, except for $\varepsilon_{\cg} = 10^{-12}$, in which case the non-flexible method cannot be distinguished from FMR.
%%%%%%%%%%%%%%%%%%%%%%%%%%%%%%%%%%%%%%%%%%%%%%%%%%%%%%

\paragraph{Acknowledgement} We are very grateful to the authors of \cite{GueLieMehSzy21}, Candan G\"ud\"uc\"u, J\"org Liesen, Volker Mehrmann and Daniel Szyld for stimulating dicussions and for making their example matrices available to us with the precise parameter settings, as well to Nicodemus Banagaaya for providing us with the matrix for Example~4 from~\cite{coupledET}.

\bibliographystyle{siamplain}
\bibliography{lit}

\end{document}